
\input epsf
\input amssym.def
\input amssym
\magnification=1100
\baselineskip = 0.25truein
\lineskiplimit = 0.01truein
\lineskip = 0.01truein
\vsize = 8.5truein
\voffset = 0.2truein
\parskip = 0.10truein
\parindent = 0.3truein
\settabs 12 \columns
\hsize = 5.4truein
\hoffset = 0.4truein

\setbox\strutbox=\hbox{%
\vrule height .708\baselineskip
depth .292\baselineskip
width 0pt}
\font\caps=cmcsc10

\font\bigtenrm=cmr10 at 14pt

\def\sqr#1#2{{\vcenter{\vbox{\hrule height.#2pt
\hbox{\vrule width.#2pt height#1pt \kern#1pt
\vrule width.#2pt}
\hrule height.#2pt}}}}
\def\square{\mathchoice\sqr46\sqr46\sqr{3.1}6\sqr{2.3}4}

\centerline{\bigtenrm NEW LOWER BOUNDS ON SUBGROUP GROWTH}
\centerline{\bigtenrm AND HOMOLOGY GROWTH}
\tenrm
\vskip 14pt
\centerline{MARC LACKENBY\footnote*{Supported by an EPSRC Advanced Fellowship}}
\vskip 18pt

\centerline{\caps 1. Introduction}
\vskip 6pt

Subgroup growth is an important new area of group theory.
It attempts to quantify the number of finite index subgroups
of a group, as a function of their index. In this paper,
we will provide new, strong lower bounds on the subgroup
growth of a variety of different groups. This will include
the fundamental groups of all finite-volume hyperbolic 3-manifolds.
By using the correspondence between
subgroups and covering spaces, we will be able to address
the following natural question: how many finite-sheeted covering 
spaces does a hyperbolic 3-manifold have, as a function of 
the covering degree? 

We will see that, when analysing the subgroup growth of a
group, it is helpful also to consider its `homology growth'.
This is concerned with the rank and order of the
first homology of its finite index subgroups. 
Fast homology growth is a useful tool when 
establishing fast subgroup growth. 

Our main result is a very general theorem, which places
a lower bound on the rank of the first homology (with mod $p$ coefficients)
of a normal subgroup $G_1$ of a group $G$, when $G/G_1$
is a finite elementary abelian $p$-group. This homology can then be
used to construct a finite index subgroup $G_2$
of $G_1$. Repeating this process, we obtain a nested
sequence of finite index subgroups $G_i$ with lower bounds
on the rank of their first homology. This works best when there is
an upper bound on the rank of the {\sl second homology}
of each $G_i$ in terms of the rank of its first homology.
Such a relationship is known to hold when $G$ is the fundamental
group of a closed 3-manifold, but it appears to be true in a much
wider context. For example, we will be able to provide
new information about the homology growth and subgroup
growth of groups with deficiency at least 1, including free-by-cyclic groups, 
and the fundamental groups of closed 4-manifolds with non-positive Euler characteristic.

An interesting aspect to this paper is that the
proofs of the main theorems are largely topological, despite
the fact that their statements are entirely algebraic
in nature.

We now give more precise statements of these results.
Let $p$ be a prime and let ${\Bbb F}_p$ be the field of
order $p$. For $r=1$ and $2$, let $b_r(G; {\Bbb F}_p)$
be the dimension of the homology group $H_r(G; {\Bbb F}_p)$. 
Thus, $b_1(G; {\Bbb F}_p)$ is the dimension over
${\Bbb F}_p$ of $G/([G,G]G^p)$, and $b_2(G; {\Bbb F}_p)$
is the mod $p$ Schur multiplier. We will be interested in groups
satisfying the following conditions.

\noindent {\bf Definition.} A group $G$ has the
{\sl $b_2-b_1$ property with respect to the prime $p$} if there is a uniform upper bound on
$$b_2(G_i; {\Bbb F}_p) - b_1(G_i; {\Bbb F}_p),$$
as $G_i$ ranges over all finite index subgroups.
A group has the {\sl $b_2/b_1$ property with respect to $p$} if there is a uniform upper
bound on $${b_2(G_i; {\Bbb F}_p) \over b_1(G_i; {\Bbb F}_p) +1 }$$
as $G_i$ ranges over all finite index subgroups.

Examples of groups satisfying the $b_2-b_1$ condition are the fundamental groups of
closed 3-manifolds and groups with deficiency at least 1. The
fundamental group of
any closed orientable 4-manifold with non-positive Euler characteristic
is a $b_2/b_1$ group. (Section 8 contains a proof of these results.)

Let $s_n(G)$ denote the number of subgroups 
of $G$ with index at most $n$.
Recall that a subgroup $K$ of a group $G$ is {\sl subnormal} 
(written $K \triangleleft \! \triangleleft \, G$) if there exists a
finite sequence of subgroups $G = G_1 \geq G_2 \geq \dots \geq G_r = K$
such that each $G_i$ is normal in $G_{i-1}$.
Let $s_n^{\triangleleft \triangleleft}(G)$
be the number of subnormal subgroups
of $G$ with index at most $n$.

The following is our main result on subgroup growth.

\noindent {\bf Theorem 1.1.} {\sl Let $G$ be a finitely generated
group that has the $b_2-b_1$ property with respect to the prime 2. Suppose that
$$\sup \{ b_1(G_i; {\Bbb F}_2) : G_i \hbox{ is a finite index subgroup of } G \} = \infty.$$
Then, for infinitely many $n$,
$$s_n(G) \geq s_n^{\triangleleft \triangleleft}(G)
> 2^{n / ( \sqrt{\log(n)} \log\log n)}.$$

}

This is a rather strong statement, since the lower
bound that it places on $s_n(G)$ and
$s_n^{\triangleleft \triangleleft}(G)$
is not far from the fastest possible
subgroup growth of a finitely generated group. It is known,
that for any finitely generated group $G$, there is a
constant $k$ such that, for all $n$,
$$\eqalign{
s_n(G) &\leq k^{n \log n} \cr
s_n^{\triangleleft \triangleleft}(G) &\leq k^n.}
$$

It will be obvious from the proof of Theorem 1.1 that the
full hypotheses of the theorem are not required. In particular,
one does not need to bound $b_2(G_i; {\Bbb F}_2) - b_1(G_i; {\Bbb F}_2)$
for all finite index subgroups $G_i$ of $G$, merely for those
in the derived 2-series of a certain finite index subgroup of $G$.
Also, one may further weaker the $b_2 - b_1$ condition, by hypothesising
that $b_2(G_i; {\Bbb F}_2) - b_1(G_i; {\Bbb F}_2)$ does not grow
too fast as a function of $b_1(G_i; {\Bbb F}_2)$. One also does
not need to assume that the supremum of $b_1(G_i; {\Bbb F}_2)$ is
infinite. For this follows from the $b_2 - b_1$ hypothesis,
provided that some $b_1(G_i; {\Bbb F}_2)$ is greater than some
constant that can be estimated. (For example, for closed 3-manifold
groups, this constant is 3.) We discuss these matters in Sections 6 and 8.

For $b_2/b_1$ groups, we can prove the following result. Although not
as strong as Theorem 1.1, it nonetheless provides good lower bounds
on subgroup growth.

Recall that for sequences $f_n$ and $g_n$, the terminology
$f_n = \Omega(g_n)$ means that $f_n/g_n \rightarrow \infty$ as
$n \rightarrow \infty$.

\noindent {\bf Theorem 1.2.} {\sl Let $G$ be a finitely generated
group that has the $b_2/b_1$ property with respect to some prime $p$. Suppose that
$$\sup \{ b_1(G_i; {\Bbb F}_p) : G_i \hbox{ is a finite index subgroup of } G \} = \infty.$$
Then, there is a constant $k > 0$ such that
$$s_n(G) \geq s_n^{\triangleleft \triangleleft}(G) = \Omega( p^{n^k}).$$}

We will find lower bounds on the constant $k$ in the above
result.

These theorems make the hypothesis that $b_1(G_i; {\Bbb F}_p)$, for suitable $p$,
can be chosen to be arbitrarily large. Clearly, some sort of hypothesis
along these lines is necessary. For example, finitely generated
abelian groups satisfy the $b_2 - b_1$ condition but have only
polynomial subgroup growth. However, this is a reasonably mild restriction,
and is often satisfied in practice.
For example, the following result is a well known consequence
of the Lubotzky alternative and the Nori-Weisfeller strong
approximation theorem (see Corollary 18 of Window 9 in [5].)

\noindent {\bf Theorem 1.3.} {\sl Let $G$ be a finitely generated linear
group. Then
either $G$ is virtually soluble, or, for any prime $p$,
$$\sup \{ b_1(G_i; {\Bbb F}_p) : G_i \hbox{ is a finite index subgroup of } G \} = \infty.$$}

In fact, it is often the case that if $b_1(G_1; {\Bbb F}_p)$ is
bigger than some fixed constant, for some finite index subgroup $G_1$ of $G$, then
the supremum of $b_1(G_i; {\Bbb F}_p)$ over all finite index
subgroups $G_i$ is infinite. To illustrate this point and to
emphasise how wide-ranging Theorems 1.1 and 1.2 are,
we give the following.

\noindent {\bf Theorem 1.4.} {\sl Let $G$ be a 
group satisfying one of the following conditions:
\item{1.} $G$ is a lattice in ${\rm PSL}(2, {\Bbb C})$;
\item{2.} $G$ is the fundamental group of a closed 3-manifold
and $b_1(G_i; {\Bbb F}_2) > 3$ for some finite index subgroup $G_i$ of
$G$;
\item{3.} $G$ has deficiency at least 1 and $b_1(G_i; {\Bbb F}_2) > 2$
for some finite index subgroup $G_i$ of $G$;
\item{4.} $G$ is a (finitely generated free non-abelian)-by-cyclic group.

\noindent Then
$$s_n(G) \geq s_n^{\triangleleft \triangleleft}(G)
> 2^{n / ( \sqrt{\log(n)} \log\log n)})$$
for infinitely many $n$.

}

\noindent {\bf Theorem 1.5.} {\sl Let $G$ be a group
satisfying one of the following conditions:
\item{1.} $G$ is the fundamental group of a closed 3-manifold
and $b_1(G_i; {\Bbb F}_p) > 3$ for some finite index subgroup $G_i$ of
$G$ and some prime $p$;
\item{2.} $G$ has deficiency at least 1 and $b_1(G_i; {\Bbb F}_p) > 2$
for some finite index subgroup $G_i$ of $G$ and some prime $p$;
\item{3.} $G$ is the fundamental group of a closed 4-manifold
with non-positive Euler characteristic and $b_1(G_i; {\Bbb F}_p) > 4$
for some prime $p$ and some finite index subgroup $G_i$ of $G$.

\noindent Then there is a constant $k > 0$ such that
$$s_n(G) \geq s_n^{\triangleleft \triangleleft}(G) = \Omega( p^{n^k}).$$}

The key piece of machinery that is the driving force behind this
paper is the following result.

\noindent {\bf Theorem 1.6.} {\sl Let $G$ be a finitely generated group, and let $p$ be
a prime such that $b_2(G; {\Bbb F}_p)$ is finite. Let $K$ be a finite index normal subgroup such that
$G/K$ is an elementary abelian $p$-group of rank $n$.
Then, for any integer $\ell$ between $0$ and $n$,
$$b_1(K; {\Bbb F}_p) \geq \sum_{r=2}^{\ell+1} \left ( n \atop r \right ) (r-1)+ 
(b_1(G; {\Bbb F}_p)-n) \sum_{r=0}^{\ell}
\left ( n \atop r \right )  - b_2(G; {\Bbb F}_p) \sum_{r=0}^{\ell-1}
\left ( n \atop r \right ).$$
Moreover, if $p=2$,
$$b_1(K; {\Bbb F}_2) \geq b_1(G; {\Bbb F}_2) \sum_{r=0}^\ell \left ( n \atop r \right ) 
- \sum_{r=1}^{\ell+1} \left ( n \atop r \right ) 
- b_2(G; {\Bbb F}_2) \sum_{r=0}^{\ell-1}
\left ( n \atop r \right ).$$
}

In the case where $\ell = 0$, these formulas should be interpreted by
taking the sums $\sum_{r=2}^{\ell+1}$ and $\sum_{r=0}^{\ell-1}$
to be zero.

Theorem 1.6 really is a collection of inequalities, one for each
integer $\ell$ between 0 and $n$, known as the `level'. In
practice, one chooses $\ell$ to obtain the strongest possible
inequality. For $p=2$, the important point here is that the
first summation runs up to $r = \ell$, whereas the third
summation goes up only to $r = \ell -1$. 

By applying this result to the derived $2$-series of some finite
index subgroup of $G$, we obtain the following lower bounds on
homology growth.

\noindent {\bf Theorem 1.7.} {\sl Let $G$ be a finitely generated group
that has the $b_2-b_1$ property with respect to the prime 2. Suppose that
$$\sup \{ b_1(G_i; {\Bbb F}_2) : G_i \hbox{ is a finite index subgroup of } G \} = \infty.$$
Then $G$ has a nested sequence of finite index normal subgroups $\{ G_i \}$,
such that
$$b_1(G_i; {\Bbb F}_2) = \Omega \left ( {[G:G_i] \over
\sqrt{ \log [G:G_i] } \log \log [G:G_i]} \right ).$$

}

This is very nearly the maximum possible growth rate of homology.
For, $b_1(G_i; {\Bbb F}_2)$ is at most the rank of $G_i$,
which, by the Reidemeister-Schreier process, is bounded
above by a linear function of $[G:G_i]$. 

The sequence $\{ G_i \}$ provided by Theorem 1.7 is the derived 2-series
of some finite index subgroup $G_1$ of $G$. Indeed, assuming that $G$ is
a finitely generated group with the $b_2-b_1$ property with respect to the prime $2$ and that $b_1(G_1; {\Bbb F}_2)$
is sufficiently large, then the derived 2-series of $G_1$ always has fast
homology growth. More precisely, for any real number $\lambda < \sqrt { 2 / \pi}$, 
there is a constant $N$, with the following property. If $G_1$ is
any finite index subgroup of $G$ with $b_1(G_1; {\Bbb F}_2) \geq N$,
and $\{ G_i \}$ is the derived 2-series of $G_1$, then
$$b_1(G_{i+1}; {\Bbb F}_2) \geq \lambda 2^{b_1(G_i;{\Bbb F}_2)} \sqrt {b_1(G_i; {\Bbb F}_2)},$$
for all $i \geq 1$. This is proved in Section 6.

Theorem 1.1 follows rapidly from Theorem 1.7, 
because if $G_i$ is a finite index normal (or just subnormal) subgroup
of $G$, then for $n = 2[G:G_i]$,
$$s_n(G) \geq s_n^{\triangleleft \triangleleft}(G) \geq 2^{b_1(G_i; {\Bbb F}_2)}.$$
Of course, Theorem 1.7 applies to any of the groups in Theorem 1.4.

A weaker form of Theorem 1.7 holds for groups satisfying
the $b_2/b_1$ condition. This applies, in particular, to any of the groups
in Theorem 1.5.

\noindent {\bf Theorem 1.8.} {\sl 
Let $G$ be a finitely generated group
that has the $b_2/b_1$ property with respect to the prime $p$. Suppose that
$$\sup \{ b_1(G_i; {\Bbb F}_p) : G_i \hbox{ is a finite index subgroup of } G \} = \infty.$$
Then $G$ has a nested sequence of finite index subgroups
$G \triangleright G_1 \triangleright G_2 \triangleright \dots$, where each
$G_i$ is normal in $G_{i-1}$, such that, for some $k > 0$,
$$b_1(G_i; {\Bbb F}_p) = \Omega({[G:G_i]^k}).$$
}

One might wonder why Theorems 1.1 and 1.7 are stated only for
the prime 2, whereas Theorems 1.2 and 1.8 work for any prime.
This is due to the asymptotics of the binomial coefficients.
A full reason is given in Section 6. Of course, however,
if a group satisfies the $b_2 - b_1$ condition with respect to
an odd prime $p$, then it also satisfies the $b_2/b_1$ condition
with respect to $p$. So, Theorems 1.2 and 1.8 provide lower
bounds on its subgroup growth and homology growth.

Theorems 1.1 and 1.7 apply to the fundamental groups of
finite-volume hyperbolic 3-manifolds. They significantly
improve the previous known lower bounds on their subgroup growth
(see Proposition 7.2.3 in [5]). They also
suggest an interesting direction for future research.
It is a major conjecture that the fundamental group of any closed hyperbolic
3-manifold should have a finite index subgroup with positive
first Betti number. Even more ambitious is the conjecture
that such groups have a finite index subgroup with a 
free non-abelian quotient. It has been shown that both
these conclusions hold for a finitely presented group $G$ provided it has an
abelian $p$-series $G \triangleright G_1 \triangleright G_2 \triangleright \dots$
with `rapid descent'. This means that each quotient $G_i/G_{i+1}$
is an elementary abelian $p$-group, and that $b_1(G_i/G_{i+1}; {\Bbb F}_p)/[G:G_i]$
is bounded away from zero. (See Theorem 1.1 of [3].)
Thus, Theorem 1.7 may represent a
first step towards a proof of these conjectures.
In any case, a good understanding of how many covering spaces
the manifold has and of their homology must surely be useful.

The paper is organised as follows. In Section 2, we establish a
preliminary technical result which produces a presentation for a group
having some useful properties. Section 3 contains the proof of Theorem 1.6,
and is the heart of the paper. Section 4 gives an explanation of
the link between Theorem 1.6, the lower central $p$-series and an
exact sequence of Stallings. In Section 5, we prove, under hypotheses
rather weaker than those in Theorems 1.1 and 1.2, that $G$ has 
finite index normal subgroups $G_i$ where $b_1(G_i; {\Bbb F}_p)$ is arbitrarily
large. In Section 6, we use the lower
bounds of Theorem 1.6 to deduce the existence of a sequence of finite
index subgroups with fast homology growth, giving Theorems 1.7 and 1.1.
In Section 7, we deal with $b_2/b_1$ groups. In Section 8,
we prove Theorems 1.4 and 1.5, which establish that our
results apply to a wide variety of different groups.

We thank the referee for carefully reading this paper, and for
suggesting several improvements to it.

\vskip 18pt
\centerline{\caps 2. Choosing a group presentation}
\vskip 6pt

Our goal in this and the next section is to prove Theorem 1.6.
Recall that we are assuming that $G$ has a normal subgroup $K$ such that
$G/K$ is an elementary abelian $p$-group of rank $n$.
In this section, we prove that the group $G$ has a presentation 
$\langle X_1, X_2, X_3 | R_1, R_2, R_3 \rangle$ having
some useful technical properties. The sets $X_1$, 
$X_2$ and $X_3$ will be finite. The free group on
the generators $X_1 \cup X_2 \cup X_3$ will be denoted by $F$.
Then $G$ is the quotient of $F$ by a normal subgroup $R$, the relations of $G$.
Let $F_-$ be the free group on the generators $X_1 \cup X_2$.

We now introduce some terminology.
For any group $H$ and integer $m \geq 1$, let $\gamma_m(H)$ be the
$m^{\rm th}$ term of the lower central $p$-series for
$H$. Recall that this is defined
recursively by setting $\gamma_1(H) = H$ and $\gamma_{m+1}(H) = 
[\gamma_m(H), H](\gamma_m(H))^p$.

Suppose that we are given an integer $m \geq 2$. (We will fix $m$ later.)
We are aiming to ensure that the presentation $\langle X_1, X_2, X_3 | R_1, R_2, R_3 \rangle$
has the following properties:
\item{(i)} $X_1$ forms a basis for $G/K$;
\item{(ii)} $X_1 \cup X_2$ forms a basis for $H_1(G; {\Bbb F}_p)$;
\item{(iii)} every element of $X_3$ is trivial in $H_1(G; {\Bbb F}_p)$;
\item{(iv)} every element of $R_2$ lies in $\gamma_2(F)$;
\item{(v)} every element of $R_3$ is of the form $x_3 = f(x_3)$,
where $f(x_3)$ is the product of an element in $\gamma_2(F_-)$ and
an element of $\gamma_m(F)$;
\item{(vi)} $R_1$ is a basis for $H_2(G; {\Bbb F}_p)$;
\item{(vii)} every element of $R_1$ lies in 
$\gamma_2(F_-) \gamma_m(F)$.

We will first prove that such a presentation can always be found.
Afterwards, we will illustrate this proof with an example.

Ensuring properties (i), (ii) and (iii) is trivial. 

We claim that, for each integer $m \geq 2$, 
$\gamma_2(F)R = \gamma_2(F_-) \gamma_m(F) R$.
This is equivalent to the statement that
$\gamma_2(G) = \gamma_2(F_-) \gamma_m(G)$. 
(This represents a slight abuse of terminology:
we are confusing $\gamma_2(F_-)$ and its image in
$G$.) Now, $G/\gamma_m(G)$ is a finite $p$-group.
The Burnside basis theorem (Theorem 12.2.1 of [2]) states that in any finite $p$-group $H$,
a set of elements generates $H$ if and only if it
generates $H/\gamma_2(H)$. Since $X_1 \cup X_2$
generates $G/\gamma_2(G)$, it therefore also generates $G/\gamma_m(G)$.
When any element of $\gamma_2(G/\gamma_m(G))$ is
expressed as a word in the generators $X_1 \cup X_2$,
its total weight in each generator is a multiple
of $p$. Thus, $\gamma_2(G) = \gamma_2(F_-) \gamma_m(G)$, as
required. This proves the claim. 

Thus, any element of $X_3$ is
equal in $G$ to the product of an element of $\gamma_2(F_-)$
and an element of $\gamma_m(F)$. When we use this fact, $m$ will be some fixed integer
at least 2, to be chosen later.
For each $x_3 \in X_3$, let $f(x_3)$ be the product
of an element of  $\gamma_2(F_-)$ and an element of $\gamma_m(F)$,
such that $x_3 = f(x_3)$ in $G$. Let $R_3$ be the relations
$\{ x_3^{-1} f(x_3) : x_3 \in X_3 \}$.

We now construct a set of relations $R_2$, as follows:
for each $r \in R$, replace each occurrence of every $x_3 \in X_3$ in $r$
by $f(x_3)$. (Note that $x_3$ may still appear in these
relations $R_2$, since $x_3$ may appear in the word
$f(x_3)$, for example.) Clearly, the subgroup of $F$ normally generated
by $R_2$ and $R_3$ is $R$.
That is, the relations $R_2 \cup R_3$ specify the
same group $G$.

We claim that, for
every relation in $R_2$, the total weight of each generator
is multiple of $p$. This is because the total weight of
every element of $X_1 \cup X_2$ in each relation is a multiple of $p$,
because these generators form a basis for $H_1(G; {\Bbb F}_p)$.
Each occurrence of an $x_3 \in X_3$ in $r \in R$ 
has been replaced by $f(x_3)$, which is an element of $\gamma_2(F)$.
In every element of $\gamma_2(F)$, the total weight of each generator
is a multiple of $p$. Thus, the claim is proved.

To summarise, we have constructed a presentation
$\langle X_1, X_2, X_3 \mid R_2, R_3 \rangle$ for $G$,
satisfying (i), (ii), (iii), (iv) and (v).

Now, the Hopf formula states that
$$H_2(G; {\Bbb F}_p) \cong {R \cap ([F,F]F^p) \over [R,F] R^p}.$$
Let $R_1'$ be a basis for $H_2(G; {\Bbb F}_p)$, where we view
each element of $R_1'$ as lying in $R \cap ([F,F]F^p)$.
From $R'_1$, we will create a new set $R_1$ of elements of
$R \cap ([F,F]F^p)$ representing the same basis for $H_2(G; {\Bbb F}_p)$,
as follows. For each occurrence of an $x_3 \in X_3$ in an element
of $R'_1$, we replace it by $f(x_3)$. Let $R_1$ be the resulting
set of relations. We claim that, when passing from $R'_1$ to $R_1$,
we have not changed the classes in $H_2(G; {\Bbb F}_p)$ that
these relations represent. This will establish that $R_1$ is
also a basis for $H_2(G; {\Bbb F}_p)$. To prove the claim,
consider one such relation $r'_1 \in R'_1$. Its total
$x_3$ weight is a multiple of $p$, since $R'_1 \subset [F,F]F^p$.
Thus, at the level of $H_2(G; {\Bbb F}_p)$, 
we have simply added a multiple of $p$ copies of
the relation $x_3^{-1}f(x_3)$ to $r'_1$. We have performed such 
an operation for each $x_3 \in X_3$. Note that $x_3^{-1} f(x_3)$ lies in $[F,F]F^p$
and therefore represents an element of $H_2(G; {\Bbb F}_p)$.
Since $H_2(G; {\Bbb F}_p)$ is an elementary abelian $p$-group, this operation therefore
does not change the class in $H_2(G; {\Bbb F}_p)$, proving the claim.

We claim that every element of $R_1$ lies in $\gamma_2(F_-)\gamma_m(F)$.
Note that in each $r'_1 \in R'_1$, the total weight of each generator
in $X_1 \cup X_2$ is a multiple of $p$. We have replaced each element
of $X_3$ with the product of an element in $\gamma_2(F_-)$ and an element
in $\gamma_m(F)$. Possibly changing each such element of $\gamma_m(F)$,
we may move it to the end of $r_1$. Thus, $r_1$ is the product of
a word in $F_-$ where each generator has weight which is a multiple of
$p$, and a word in $\gamma_m(F)$. Thus, $R_1 \subset \gamma_2(F_-) \gamma_m(F)$,
as required.

Adding in
this set of relations $R_1$ to the presentation above, we obtain
a presentation
$$G = \langle X_1, X_2, X_3 \mid R_1, R_2, R_3 \rangle,$$
satisfying (i) - (vii) above, as required.

We now give an example that demonstrates the construction of such a presentation
in practice. Let $G$ be
$$\langle x_1, x_2, x_3 \mid x_3^{-1} [x_3,x_1], x_1^2 x_3^2\rangle.$$
Then $H_1(G; {\Bbb F}_2)$ has rank 2, with $x_1$ and $x_2$ as
generators. Let $K$ be the subgroup of $G$ generated by $x_2$ and 
$[G,G]G^2$. Thus, $G/K$ is isomorphic to ${\Bbb Z}/2{\Bbb Z}$.
Let $m$ be $3$, say.

We start by setting $X_1 = \{ x_1 \}$, $X_2 = \{ x_2 \}$ and $X_3 = \{ x_3 \}$.
Then, clearly, (i), (ii) and (iii) are satisfied.
In the next stage of the procedure, we find $f(x_3)$, which
equals $x_3$ in $G$ and which is the product of an element of $\gamma_2(F_-)$ and an element of $\gamma_3(F)$.
Since $x_3 = [x_3,x_1]$ (by the first relation in the group),
we may insert $x_3$ into this commutator to obtain the
relation $x_3 = [[x_3, x_1], x_1]$. We set $f(x_3)$ to be
$[[x_3, x_1], x_1]$, which is the product of an element of 
$\gamma_2(F_-)$ (the identity) and an element of $\gamma_3(F)$, as required.
Thus,
$$R_3 = \{ x_3^{-1} [[x_3, x_1], x_1] \}.$$
To construct $R_2$, we need to substitute every occurrence
of $x_3$ in the relations with $f(x_3)$. It suffices to do this
for the initial defining relations for the group.
Thus, we may set
$$R_2 = \{ [[x_3, x_1], x_1]^{-1} [[[x_3, x_1], x_1], x_1], x_1^2 [[x_3, x_1], x_1]^2 \}.$$
For the final stage of the procedure, we start with a basis for $H_2(G; {\Bbb F}_2)$.
In this case, $H_2(G; {\Bbb F}_2) \cong {\Bbb F}_2$ with generator $x_1^2 x_3^2 \in R \cap ([F,F]F^2)$.
Set $R'_1$ to be this generator. To create $R_1$ from $R'_1$, we substitute every occurrence
of $x_3$ with $f(x_3)$. Thus,
$$R_1 = \{ x_1^2 [[x_3, x_1], x_1]^2 \}.$$

We now explain how the integer $m$ is chosen.

\noindent {\bf Lemma 2.1.} {\sl For some integer $m \geq 2$,
$\gamma_m(F) \subset \gamma_2(K)$.}

\noindent {\sl Proof.} Note that $\gamma_2(K)$ is a normal
subgroup of $F$. Its index is a power of $p$, and so $F/\gamma_2(K)$
is a finite $p$-group. The lower central $p$-series of any finite $p$-group
terminates. Hence, for all sufficiently large integers $m$,
$\gamma_m(F) \subset \gamma_2(K)$. $\square$

We now fix $m$ to be the integer given by the above lemma.
(We can take $m$ to be ${\rm rank}(G/K)+2$, but we will not need this fact.)

Define subgroups of $R$ recursively, by setting
$R_{(1)} = R$ and letting $R_{(j+1)} = [R_{(j)},F](R_{(j)})^p$,
for each $j \geq 1$. An elementary induction establishes that
each $R_{(j)}$ is normal in $F$.

\noindent {\bf Lemma 2.2.} {\sl 
Let $S = \langle \! \langle R_1, R_3 \rangle \! \rangle$, the subgroup
of $F$ normally generated by $R_1$ and $R_3$. Then,
for each $j \geq 1$, $R = S R_{(j)}$.}

\noindent {\sl Proof.} We first establish the
inclusions
$$\eqalign{
R &= \langle \! \langle R_1 \rangle \! \rangle
\langle \! \langle R_2 \rangle \! \rangle
\langle \! \langle R_3 \rangle \! \rangle \cr
&\subseteq (R \cap ([F,F]F^p)) 
\langle \! \langle R_3 \rangle \! \rangle \cr
&= \langle R_1 \rangle [R,F]R^p
\langle \! \langle R_3 \rangle \! \rangle \cr
&\subseteq S [R,F]R^p \cr
&\subseteq R.}$$
In the second step above, we use are using properties
(iv) and (vii) of the presentation. In the third
step, we are using property (vi). We deduce that
each of these inclusions is an equality.

We now prove the lemma by induction on $j$. For $j =1$,
this is trivial. Suppose that it is true for a given $j$. Then
$$\eqalign{
R &= S [R,F]R^p \cr
&= S[SR_{(j)}, F] (SR_{(j)})^p \cr}$$
$$\eqalign{
&= S[S,F][R_{(j)},F] S^p (R_{(j)})^p \cr
&= S[R_{(j)},F](R_{(j)})^p \cr
&= SR_{(j+1)}.}$$
In the first equality, we are using the claim proved above. In the
second, the inductive hypothesis is used. In the third equality,
we are using the fact that $[AB,C] = [A,C][B,C]$, for normal
subgroups $A$, $B$ and $C$ of a group.
$\square$

\vskip 18pt
\centerline {\caps 3. Finding homology in covering spaces}
\vskip 6pt

The aim of this section is to prove Theorem 1.6, which provides
lower bounds on $b_1(K; {\Bbb F}_p)$, for certain subgroups $K$
of a group $G$. We work with a presentation $\langle X_1, X_2, X_3 \mid
R_1, R_2, R_3 \rangle$ for $G$, satisfying conditions (i) - (vii)
of Section 2. This determines a 2-complex $L$, in the usual way.
It has a single 0-cell, which we take to be its basepoint, an
oriented 1-cell for each generator, and a 2-cell for each relation.
Then $\pi_1(L)$ is isomorphic to $G$. We are considering a normal
subgroup $K$ of $G$, such that $G/K$ is an elementary abelian
$p$-group of rank $n$, for some prime $p$. Associated with $K$, there is a finite-sheeted
cover $\tilde L$ of $L$. We fix a basepoint for $\tilde L$ that maps to
the 0-cell of $L$. The plan is to find a lower bound
on $b_1(\tilde L; {\Bbb F}_p)$. Since $b_1(\tilde L; {\Bbb F}_p)$ equals $b_1(K; {\Bbb F}_p)$,
this will give the required lower bound on $b_1(K; {\Bbb F}_p)$. Now,
$b_1(\tilde L; {\Bbb F}_p)$ is the dimension of $H^1(\tilde L; {\Bbb F}_p)$.
Recall that this is the vector space $Z^1/B^1$. Here,
$B^1$ and $Z^1$ are subspaces of $C^1$, the space of
all 1-cochains on $\tilde L$, with mod $p$ coefficients.
The 1-coboundaries $B^1$ are precisely those 1-cochains
such that their evaluation on any closed loop is trivial.
The 1-cocycles $Z^1$ are those 1-cochains with trivial
evaluation on the boundary of any 2-cell. Thus, $Z^1$ is a subspace
of $C^1$, obtained by imposing $p^n |R_1 \cup R_2 \cup R_3|$ linear
constraints, since $p^n |R_1 \cup R_2 \cup R_3|$ is the
number of 2-cells of $\tilde L$ (which may be infinite).
In our approach to the proof of Theorem 1.6, we work with
certain subspaces of $C^1$, which we denote by $U^1_\ell$,
for integers $\ell$ between 0 and $n$. (The integer $\ell$ is
the same as that in the statement of Theorem 1.6, and is known
as the `level' of the subspace $U^1_\ell$.) These subspaces form
a hierarchy
$$C^1 \supseteq U^1_n \supseteq U^1_{n-1} \supseteq \dots
\supseteq U^1_0.$$
Each has the following dimension:
$${\rm dim}(U^1_\ell) = b_1(G; {\Bbb F}_p) \sum_{i=0}^\ell \left ( n \atop i \right ).$$
These subspaces have the following nice property. Roughly
speaking, to determine
whether or not certain cochains in $U^1_\ell$ are in fact cocycles, one
does not need to verify that their evaluation around every 2-cell
is zero. Instead, it suffices to check a certain set of
2-cells, with cardinality at most
$$b_2(G; {\Bbb F}_p) \sum_{i=0}^{\ell-1} \left ( n \atop i \right ).$$
More precisely, if an element of $U^1_\ell$ has zero
evaluation around these 2-cells, then a cocycle may
be constructed from it.
Of course, some of these cocycles may be coboundaries, but in fact 
this is true only for a relatively small subspace.

For any $j \in X_1 \cup X_2$, let $c_j$ be
the 1-cochain on $L$ which sends the edge
of $L$ labelled $j$ to $1 \in {\Bbb F}_p$, and maps the remaining
edges to 0. These are cocycles because $X_1 \cup X_2$
forms a basis for $H_1(G; {\Bbb F}_p)$. 

We lift the orientations on the 1-cells of $L$ to the 1-cells of
$\tilde L$. For each 1-cell $e$ of $\tilde L$, let $i(e)$ denote its
initial vertex. The map $\tilde L \rightarrow L$ sends $e$ to a 1-cell of $L$, and so
$e$ has a well-defined evaluation under each $c_j$. We denote this by
$c_j(e)$. 

Every vertex $v$ of $\tilde L$ also has a well-defined evaluation
under $c_j$, for each $j \in X_1$, defined as follows. Pick a path
from the basepoint of $\tilde L$ to $v$. This projects to a loop
$\alpha$ in $L$. Define $c_j(v)$ to be $c_j(\alpha)$. This is
well-defined, because if $\alpha'$ is another path from the basepoint
to $v$, then
$$c_j(\alpha') = c_j(\alpha'. \alpha^{-1}) + c_j(\alpha) = c_j(\alpha).$$
The final equality holds because $\alpha'.\alpha^{-1}$ is a loop
in $\tilde L$ and so $\alpha'. \alpha^{-1} \in K$, which implies
that $c_j(\alpha'.\alpha^{-1}) = 0$.

We can now define the subspace $U^1_\ell$ of $C^1$, for each integer $\ell$
between 0 and $n$. We do this by specifying a spanning set. For each
subset $A$ of $X_1$ with size at most $\ell$, and for
each element $y \in X_1 \cup X_2$, define $c(A,y)$
to be the following 1-cochain. On an edge $e$ of $\tilde L$, let
$$c(A,y)(e) = \left ( \prod_{j \in A} c_j(i(e)) \right ) c_{y}(e).$$
When $A = \emptyset$, we take this to mean that
$c(A,y)(e) = c_y(e)$, by convention.
Then $U^1_\ell$ is defined to be the subspace of $C^1$ spanned by
these cochains.

\noindent {\bf Example.} Let $G$ be the free group on 3 generators,
let $K = [G,G]G^2$ and let $p=2$. Then $L$ is the wedge of 3 circles, and
$\tilde L$ is the cube-shaped graph shown twice in Figure 1. Note that
the dotted edges in each figure join up with each other. Then
the support of the cochains $c(\{ 1,2 \}, 3)$ and $c(\{1 \}, 1)$ is 
shown in bold.

\vskip 18pt
\centerline{\epsfbox{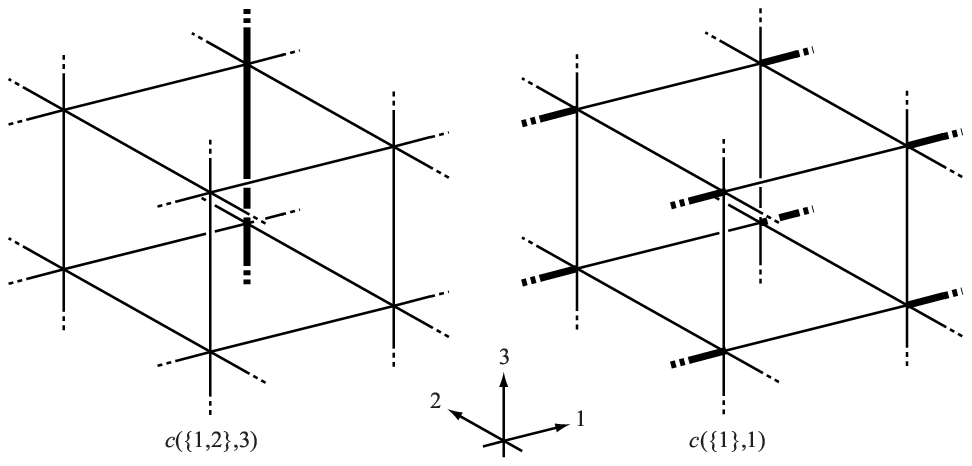}}
\vskip 18pt
\centerline{Figure 1.}

It is possible to prove that the cochains $c(A,y)$ are linearly 
independent and so form a basis for $U^1_\ell$.
This is a reassuring result. But since we will not use this fact,
we omit its proof.

Given any element $g \in F$ and a 1-cochain $c$ on $\tilde L$, we define $c(g)$ as follows. Pick a word
representing $g$, which specifies a path in $\tilde L$ starting at
the basepoint. Define $c(g)$ to be the evaluation of $c$
on this path. This is clearly independent of the choice of word
representing $g$.

Pick a total ordering on $X_1$.
For a subset $E$ of $X_1$, with elements
$i_1 < i_2 < \dots < i_{|E|}$, let $w_E$ be the word
$i_1 \dots i_{|E|}$ representing an element of the free group $F$. When $E = \emptyset$,
then $w_E$ is the identity element of $F$. 

Our key technical result is the following.

\noindent {\bf Proposition 3.1.} {\sl Let $z$ be an element
of $U^1_\ell$. Let $k$ be an element of $R$. Suppose that $z(w_E k w_E^{-1}) = 0$,
for all subsets $E$ of $X_1$ with size at
most $\ell -1$. Then $z(gkg^{-1}) = 0$
for all $g \in F$.}

We prove this using a series of lemmas. Throughout,
$A$ is a subset of $X_1$ with size at most $\ell$, and $y$
is an element of $X_1 \cup X_2$. We define $KR$ to be the subgroup of $F$
that maps to $K$ under the quotient homomorphism $F \rightarrow G$.
Note that $KR$ is the set of elements of $F$ that form closed loops in $\tilde L$.

\noindent {\bf Lemma 3.2.} {\sl For $k \in KR$ and $g \in F$,
$$c(A,y)([g,k]) = 
\sum_{\scriptstyle B \subseteq A \atop \scriptstyle B\not= \emptyset} 
\left ( c(A - B,y)(k) \prod_{j \in B} c_j(g) \right ).$$
}

\noindent {\sl Proof.} 
Represent $k$ by the word $x_1 \dots x_s$ in the
generators of $F$. Because this is a closed loop in $\tilde L$,
the parts of $[g,k]$ in $g$ and $g^{-1}$ run along the
same edges in reverse and hence cancel. Thus, we need only consider 
evaluation of $c(A,y)$ on the parts of $[g,k]$ in $k$ and $k^{-1}$. 
Let $e_1, \dots, e_s$ be the edges of $k$. Then the $k$
part of $[g,k]$ runs along $ge_1, \dots, ge_s$, where each
$ge_r$ denotes the copy of $e_r$ translated by the covering
transformation of $\tilde L$ corresponding to $g$.
The $k^{-1}$ part of $[g,k]$ runs along $e_1, \dots, e_s$ in reverse.
Let $\beta_r$ be $\pm 1$, according to whether $k$ runs
forwards or backwards along the edge $e_r$.
So,
$$\eqalign{
& \qquad c(A,y)([g,k]) \cr
&=  \sum_{r=1}^{s} \beta_r c(A,y) (g e_r) 
- \sum_{r=1}^{s} \beta_r c(A,y)(e_r)
\cr
& = \sum_{r=1}^{s} \left ( 
\left ( \prod_{j \in A}  c_j(i(g e_r)) \right ) \beta_r c_y(e_r)
- \left ( \prod_{j \in A} 
c_j(i(e_r)) \right ) \beta_r c_y(e_r)
 \right) \cr
& = \sum_{r=1}^{s} \left ( 
\prod_{j \in A} 
(c_j(g) + c_j(i(e_r))) 
- \prod_{j \in A} 
c_j(i(e_r)) \right ) \beta_r c_y(e_r) \cr
&= \sum_{r=1}^{s} \Biggl (
\sum_{B \subseteq A} 
\Big ( \prod_{j \in B} c_j(g) \prod_{j \in A - B}
c_j(i(e_r)) \Big) - \prod_{j \in A} 
c_j(i(e_r)) \Biggr ) \beta_r c_y(e_r)\cr
&= \sum_{\scriptstyle B \subseteq A \atop \scriptstyle B \not= \emptyset}
\Biggl ( \prod_{j \in B} c_j(g) \Biggr )
\Biggl ( \sum_{r=1}^{s}
\Big ( \prod_{j \in A - B} c_j(i(e_r)) \Big ) \beta_r c_y(e_r)
\Biggr ) \cr
&= \sum_{\scriptstyle B \subseteq A \atop \scriptstyle B \not= \emptyset}
\left ( c(A - B,y)(k) \prod_{j \in B} c_j(g) \right ),}$$
as required. $\square$

\noindent {\bf Corollary 3.3.} {\sl For any $k \in KR$
and $j \in X_1 \cup X_2 \cup X_3$,
$$c(A,y)([j,k]) = \cases{
c(A - \{ j \},y)(k) & if $j \in A$, \cr
0 & otherwise.}$$}

\noindent {\bf Lemma 3.4.} {\sl For any $E \subseteq X_1$, and any $k \in KR$,
$$c(A,y)(w_E k w_E^{-1}) = \sum_{B \subseteq A \cap E}
c(A - B,y)(k).$$}

\noindent {\sl Proof.}
We prove this by induction on $|E|$. For $E = \emptyset$,
it is trivial. For the inductive step,
let $j$ be the first element of $E$, and let $D$ be
$E - \{ j \}$.
Then 
$$c(A,y)(w_E k w_E^{-1}) 
= c(A,y)(j w_D k w_D^{-1} j^{-1}).$$
If $j$ is not in $A$, then, by Corollary 3.3,
this equals
$$c(A,y)(w_D k w_D^{-1}) = 
\sum_{B \subseteq A \cap D}
c(A - B,y)(k)
= \sum_{B \subseteq A \cap E}
c(A - B,y)(k),$$
as required. On the other hand, if
$j$ is in $A$, then it equals
$$\eqalign{
&\quad c(A,y)(w_D k w_D^{-1}) + c(A,y)([j,w_D k w_D^{-1}]) \cr
&= c(A,y)(w_D k w_D^{-1}) + c(A - \{j \},y)(w_D k w_D^{-1}) \cr
&= \sum_{B \subseteq A \cap D}
c(A - B,y)(k) + \sum_{B \subseteq (A - \{ j \}) \cap D}
c(A - \{ j \} - B,y)(k) \cr
& = \sum_{B \subseteq A\cap E}
c(A - B,y)(k).}$$
The first equality above relies on Corollary 3.3.
The second equality uses the inductive hypothesis. The induction is
established. $\square$

\noindent {\bf Lemma 3.5.} {\sl Let $z$ be an element of $U_\ell^1$. Thus,
$$z = \sum_{A,y} \lambda_{A,y} c(A,y),$$
where $A$ ranges over all subsets of $X_1$
with size at most $\ell$ and $y \in X_1 \cup X_2$,
and where $\lambda_{A,y}$ are coefficients in 
${\Bbb F}_p$. Let $k$ be an element of $R$. Suppose that 
$z(w_E k w_E^{-1}) = 0$ for all subsets $E$ of
$X_1$ with size at most $\ell - 1$. Then, for any subset 
$E$ of $X_1$,
$$\sum_{\scriptstyle \{ A: E \subseteq A, \ |A| \leq \ell \}
\atop \scriptstyle y} \lambda_{A,y} c(A - E,y)(k) = 0.$$
}

\noindent {\sl Proof.} 
We prove the lemma by induction on $|E|$. For $|E| = 0$,
it is trivial. The inductive step, when $|E| \leq \ell-1$,
follows from the fact that
$$\eqalign{
0 &= z(w_E k w_E^{-1}) \cr
&= \sum_{A,y} \lambda_{A,y} c(A,y)(w_E k w_E^{-1}) \cr
& = \sum_{A,y} \lambda_{A,y} \sum_{B \subseteq A \cap E}
c(A - B,y)(k)\cr
& = \sum_{B \subseteq E} \sum_{\scriptstyle \{ A: B \subseteq A, \ |A| \leq \ell \}
\atop \scriptstyle y}
\lambda_{A,y} c(A - B,y)(k). \cr}$$
The third equality is an application of Lemma 3.4.
By induction
$$\sum_{\scriptstyle \{ A: B \subseteq A, \ |A| \leq \ell \}
\atop \scriptstyle y}
\lambda_{A,y} c(A - B,y)(k)$$
is zero, when $B$ is strictly contained
in $E$. Hence,
$$\sum_{\scriptstyle \{ A: E \subseteq A, \ |A| \leq \ell \}
\atop \scriptstyle y}
\lambda_{A,y} c(A - E,y)(k)$$
is also zero, as required.
For $|E| \geq \ell$, note that if $E \subseteq A$, then
$E = A$. Thus, in this case, the formula we must prove is
$\sum_y \lambda_{E,y} c(\emptyset,y)(k) = 0$. But, $ c(\emptyset,y)(k) = c_y(k)$,
which is zero because we are assuming that $k \in R$.
$\square$

\noindent {\bf Lemma 3.6.} {\sl Let $z$ and $k$ be as in Lemma 3.5. Thus,
$$z = \sum_{A,y} \lambda_{A,y} c(A,y),$$
where $A$ ranges over all subsets of $X_1$
with size at most $\ell$ and $y \in X_1 \cup X_2$,
and where $\lambda_{A,y}$ are coefficients in 
${\Bbb F}_p$.
Then for any $g \in F$ and any $E \subseteq X_1$,
$$\sum_{\scriptstyle \{ A: E \subseteq A, \ |A| \leq \ell \}
\atop \scriptstyle y} \lambda_{A,y} c(A - E,y)(gkg^{-1}) = 0.$$}

\noindent {\sl Proof.}
We represent $g$ by a word in the generators.
We prove the lemma by induction on the length of this word. We already
know it to be the case when $g$ is the identity, by Lemma 3.5.
For the inductive step, applied to some word $g = j u$:
$$\eqalign{
& \qquad \sum_{\scriptstyle \{ A: E \subseteq A, \ |A| \leq \ell \}
\atop \scriptstyle y} \lambda_{A,y} c(A - E,y)(gkg^{-1}) \cr
&= \sum_{\scriptstyle \{ A: E \subseteq A, \ |A| \leq \ell \}
\atop \scriptstyle y} \lambda_{A,y} c(A - E,y)(ju ku^{-1} j^{-1}) \cr
&= \sum_{\scriptstyle \{ A: E \subseteq A, \ |A| \leq \ell \}
\atop \scriptstyle y} \lambda_{A,y} c(A - E,y)(u ku^{-1}) \cr
&+ \sum_{\scriptstyle \{ A: E \subseteq A, \ |A| \leq \ell \}
\atop \scriptstyle y} \lambda_{A,y} c(A - E,y)([j, u ku^{-1}]) \cr
&= \sum_{\scriptstyle \{ A: E \subseteq A, \ |A| \leq \ell \}
\atop \scriptstyle y}
\lambda_{A,y} c(A - E,y)(u ku^{-1}) \cr
&+ \sum_{\scriptstyle \{ A: E \subseteq A, \ |A| \leq \ell, \ j \in A - E  \}
\atop \scriptstyle y}
\lambda_{A,y} c((A - E) - \{ j \})(u ku^{-1}), \cr}$$
by Corollary 3.3.
The first sum is zero by induction, as is the second, since it
equals
$$\sum_{\scriptstyle \{ A: E \cup \{ j \} \subseteq A, \ |A| \leq \ell \}
\atop \scriptstyle y} \lambda_{A,y} c(A - (E \cup \{j \}))(u ku^{-1}).$$
$\square$

\noindent {\sl Proof of Proposition 3.1.} Set $E = \emptyset$ in
Lemma 3.6. We obtain
$$z(gkg^{-1}) = \sum_{A,y} \lambda_{A,y} c(A,y)(gkg^{-1}) = 0$$
for all $g \in F$. $\square$

So far, we have focused on cochains supported on edges
labelled by generators in $X_1 \cup X_2$. Let $U^1$ be
the space of all such cochains. We now show
how a cochain in $U^1$ has a natural modification, which has
the same values on the edges labelled by $X_1 \cup X_2$, but
so that its support might also include edges
labelled by $X_3$. This modification depends on the presentation 
for $G$ that we fixed in Section 2. We will define 
a linear map $\psi \colon U^1 \rightarrow C^1$. The modification of
a cochain $z$ will be $\psi(z)$. We define $\psi(z)$ to agree with $z$
on the edges of $\tilde L$ labelled by $X_1 \cup X_2$. Each remaining edge
$e$ is labelled by an element $x_3 \in X_3$. For this
$x_3$, there is a relation of the form $x_3 = f(x_3)$
in $R_3$. Note that $e$ is a loop based at the vertex $i(e)$. Let $g$ be a word in $X_1$
which specifies a path from the basepoint of $\tilde L$
to $i(e)$. We define $(\psi(z))(e)$ to be $z(gf(x_3)g^{-1})$.
This is clearly independent of the choice of 
$g$, since the $g$ and $g^{-1}$ parts of the loop $gf(x_3)g^{-1}$ 
traverse the same edges in the opposite direction. Note that $\psi$ is an injection.

\noindent {\bf Lemma 3.7.} {\sl The restriction of 
$\psi$ to $B^1$ is the identity. Hence, $\psi(B^1) = B^1$.}

\noindent {\sl Proof.} Note first that each element of
$B^1$ is supported on edges labelled by $X_1 \cup X_2$,
and so $B^1$ lies in $U^1$, which is the domain of $\psi$. Thus,
it makes sense to speak of the restriction of
$\psi$ to $B^1$.

Let $z$ be an element of $B^1$. On each edge $e$ labelled by
$x_3 \in X_3$, $(\psi(z))(e)$ is defined to be the
evaluation of $z$ on $gf(x_3)g^{-1}$, for suitable
$g \in F$. This is a closed loop, and so its evaluation
under the coboundary $z$ is trivial. Thus, $\psi(z)$
is trivial on all edges labelled by $X_3$. On the remaining
edges, $\psi(z)$ and $z$ agree. Thus,
$\psi(z) = z$, as required. $\square$

We now focus on a subspace of $U^1_\ell$. This will arise
as the kernel of a linear map $\phi_\ell \colon U^1_\ell \rightarrow V_\ell$,
where $V_\ell$ is the product of
$$b_2(G; {\Bbb F}_p) \sum_{i=0}^{\ell-1}
\left ( n \atop i \right )$$
copies of ${\Bbb F}_p$. We parametrise the co-ordinates
of each element of $V_\ell$ by pairs $(r_1, E)$, where $r_1 \in R_1$ and
$E \subset X_1$, with $|E| \leq \ell-1$.
The map $\phi_\ell$ is defined as follows.
If $z \in U^1_\ell$, then for each pair $(r_1, E)$, 
the $(r_1, E)$ co-ordinate of $\phi_\ell(z)$ is
 $z(w_E r_1 w_E^{-1})$.

We now define the subspace $C^1_\ell$ of $C^1$ to
be $\psi({\rm ker}(\phi_\ell))$. Thus, to construct an
element of $C^1_\ell$, start with a linear combination
$z$ of the cochains $c(A,y)$. Evaluate $z$ on
the words $w_E r_1 w_E^{-1}$. Restrict attention
to those $z$ that have zero evaluation on these words.
Given such a $z$, modify it to $\psi(z)$, which
assigns certain values to edges labelled by $X_3$.
Each such $\psi(z)$ is an element of $C^1_\ell$,
and conversely each element of $C^1_\ell$ is constructed
in this way.

\noindent {\bf Theorem 3.8.} {\sl Each element of
$C^1_\ell$ is a cocycle.}

Before we prove this, we need a lemma. Recall that $m$ is the integer
from Lemma 2.1.

\noindent {\bf Lemma 3.9.} {\sl For all elements
$g \in \gamma_m(F)$ and all cochains $c$ in $C^1$,
$c(g) = 0$.}

\noindent {\sl Proof.} This is simply a restatement of
Lemma 2.1 in topological language. $\square$

\noindent {\sl Proof of Theorem 3.8.}  Let $z$ be an element
of ${\rm ker}(\phi_\ell)$. We want to prove that
$\psi(z)$ is a cocycle. Lemma 2.2 gives that
$R = SR_{(m)}$, where $S = \langle \! \langle R_1, R_3 \rangle \! \rangle$.
Hence, to show that $(\psi(z))(r) = 0$ for all $r \in R$,
it suffices to check that $(\psi(z))(r) = 0$ for all $r \in R_{(m)} \cup S$.
But, $R_{(m)}$ lies in $\gamma_m(F)$, by the definition of
$R_{(m)}$. Lemma 3.9 implies that the evaluation of any
1-cochain in $\tilde L$ on an element of $\gamma_m(F)$ is trivial.
Thus, $(\psi(z))(r) = 0$ for all $r \in R_{(m)}$.
We therefore only need to prove that $(\psi(z))(r) = 0$
for all $r \in \langle \! \langle R_1 \rangle \! \rangle
\cup \langle \! \langle R_3 \rangle \! \rangle$.

\noindent {\sl Claim 1.} Let $e$ be an edge labelled by $x_3 \in X_3$,
and let $g$ be as in the definition of
$(\psi(z))(e)$. Then,
$(\psi(z))(gf(x_3)g^{-1}) = z(gf(x_3)g^{-1})$. 

Now, $f(x_3)$ is the product of an element of $\gamma_2(F_-)$
and an element of $\gamma_m(F)$. The evaluation of the latter
under any 1-cochain in $\tilde L$ is trivial, by Lemma 3.9.
Thus, $(\psi(z))(gy_3g^{-1})$ is equal to the evaluation under
$\psi(z)$ of the word in $\gamma_2(F_-)$, conjugated by $g$.
This lies in $F_-$, and hence the corresponding loop is supported
on the edges labelled by $X_1 \cup X_2$. But $z$ and $\psi(z)$ agree on these
edges. This proves the claim.

The claim implies that, for each 2-cell of $\tilde L$ labelled by
an element of $R_3$, the evaluation of $\psi(z)$ around its
boundary is zero. Hence, $(\psi(z))(r) = 0$ for all
$r \in \langle \! \langle R_3 \rangle \! \rangle$.

\noindent {\sl Claim 2.} For each 2-cell of $\tilde L$ labelled by an
element in $R_1$, the evaluation of $\psi(z)$ around its
boundary is zero.

By property (vii) of Section 2, each such element lies in
$\gamma_2(F_-) \gamma_m(F)$. Thus, by the argument of Claim 1,
its evaluation under $\psi(z)$ equals its evaluation under
$z$. Now, $z$ lies in ${\rm ker}(\phi_\ell)$ and so its evaluation
on each word $w_E r_1 w_E^{-1}$ is zero (where $r_1 \in R_1$ and
$E \subseteq X_1$ with $|E| \leq \ell - 1$). Proposition
3.1 then implies that its evaluation on $gr_1g^{-1}$ for any
$r_1 \in R_1$ and $g \in F$ is zero. This proves the claim
and the theorem.
$\square$

Theorem 3.8 establishes that the cochains in $C^1_\ell$ are
cocycles. But in order to prove Theorem 1.6, we need to know
how many of these cocycles are coboundaries. We start by
examining which elements of $U^1_\ell$ are coboundaries.
There will, in general, be some non-zero coboundaries, as the following
example demonstrates.

\vskip 18pt
\centerline{
\epsfxsize=3in
\epsfbox{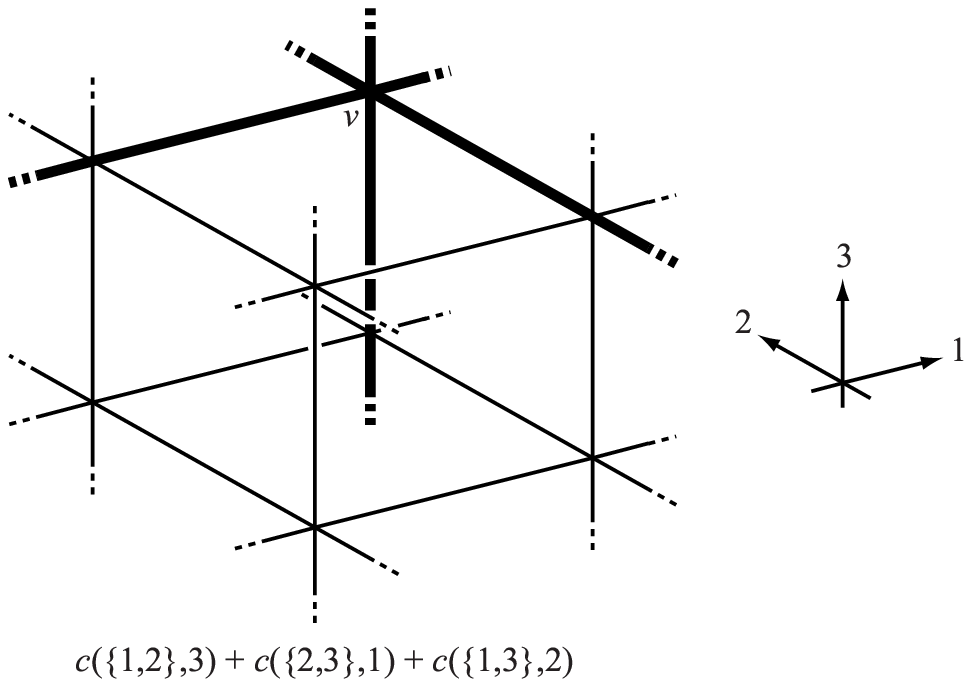}}
\vskip 18pt
\centerline{Figure 2.}
\vfill\eject

\noindent {\bf Example.} Let $G$ be the free group on 3 generators,
let $K$ be $[G,G]G^2$, and let $p=2$. Then, on $\tilde L$,
$$c(\{1,2 \}, 3) + c(\{ 2,3 \},1) + c( \{1,3 \},2)$$
is a coboundary (see Figure 2). It is the coboundary of the function supported at the vertex
$v$ where $c_1(v) = c_2(v) = c_3(v) = 1$.

In fact, it is not hard to show that, more generally, if $A$
is any non-empty subset of $X_1$, then
$$\sum_{y \in A} c(A - \{y \}, y)$$
is a coboundary, although we will not need this fact.

\noindent {\bf Proposition 3.10.} {\sl The dimension of
$U^1_\ell/(B^1 \cap U^1_\ell)$ is at least
$$\sum_{r=2}^{\ell+1} \left ( n \atop r \right )(r-1) + 
(b_1(G; {\Bbb F}_p)-n) \sum_{r=0}^{\ell}
\left ( n \atop r \right )$$
if $p$ is odd, and at least
$$b_1(G; {\Bbb F}_2) \sum_{r=0}^\ell \left ( n \atop r \right ) 
- \sum_{r=1}^{\ell+1} \left ( {n \atop r} \right )$$
if $p=2$.}

We prove this by evaluating elements of $U^1_\ell$ on a
certain set $T$ of `test' loops in $\tilde L$. 
The number of such loops will
be equal to the quantities given in Proposition 3.10. 
Let $U^1_\ell \rightarrow {\Bbb F}_p^T$ be the map that
sends an element of $U^1_\ell$ to its evaluation under the
test loops $T$. We will show that this map has rank $|T|$.
Since this map factors through $U^1_\ell/(B^1 \cap U^1_\ell)$,
this will prove the proposition.

Let $A$ be a non-empty subset of $X_1$ with cardinality at most $\ell$.
Let $y$ be an element of $(X_1 \cup X_2) - A$.
Let $y_1$ be the smallest element of $A$.
(Recall that $X_1$ comes with a total ordering.) 
We insist that if $y$ lies in $X_1$, then it is larger
than $y_1$. Define $t(A,y)$ to be the loop $[y, y_1^{-1}]$
which starts at the vertex $v$ such that
$$c_j(v) = \cases{ 1 & if $j \in A$; \cr 0 & if $j \not \in A$.}$$
When $y \in X_2$, define $t(\emptyset, y)$ to be the loop $y$
which starts and ends at the basepoint of $\tilde L$.
When $p$ is odd, the set $T$ of test loops will be all such $t(A, y)$. 
When $p =2$, the test loops will be all these loops, together with 
the following. Let $A$ be as above, but now let $y$ be
an element of $A$. Define $t(A,y)$ to be the loop $y^2$ which starts 
at the vertex $v$ defined above.
We include all such $t(A,y)$ as test loops.

\noindent {\bf Lemma 3.11.} {\sl The number of test loops $T$
is given by the quantities in Proposition 3.10.}

\noindent {\sl Proof.} Let us first consider the case
where $p$ is odd. We wish to count the number of pairs
$(A,y)$ satisfying the above conditions. 
If $y \in X_1$, then we simply count the possibilities
for $A \cup y$. For each such set $A \cup y$, with
cardinality $r$ between $2$ and $\ell+1$, there are $r-1$
choices for $y$, since $y$ cannot be the smallest
element of $A \cup y$. Thus, the number of such pairs
$(A,y)$ where $y \in X_1$ is
$$\sum_{r=2}^{\ell+1} \left ( n \atop r \right )(r-1).$$
The number of pairs $(A,y)$ where $y \in X_2$ is clearly
$$(b_1(G; {\Bbb F}_p)-n) \sum_{r=0}^{\ell}
\left ( n \atop r \right )$$
since $|X_2| = (b_1(G; {\Bbb F}_p)-n)$.

Let us now examine the case where $p=2$. Here,
we count all pairs $(A,y)$ where $A \subset X_1$,
$|A| \leq \ell$ and $y \in X_1 \cup X_2$. We then
subtract off the number that do not satisfy the
given condition. The first count gives
$$b_1(G; {\Bbb F}_2) \sum_{r=0}^\ell \left ( n \atop r \right ).$$
If $(A,y)$ does not satisfy the condition
required to define a test loop, then $y$ lies in
$X_1$ and is strictly smaller than every element of $A$. Thus,
given $A \cup y$, it is possible to determine $y$.
Thus, we need only count the number of possibilities
for $A \cup y$:
$$\sum_{r=1}^{\ell+1} \left ( n \atop r \right ).$$
The required formula follows immediately. $\square$

\vfill\eject
\noindent {\bf Lemma 3.12.} {\sl Let $A'$ be a subset of $X_1$ with
size at most $\ell$, and let $y'$ be an element of $X_1 \cup X_2$.
Let $A$ be a subset of $X_1$ with size at most $|A'|$,
and let $y$ be an element of $(X_1 \cup X_2) - A$.
Suppose that, if $y$ lies in $X_1$, then $A$ is non-empty and
$y$ is larger than the smallest element $y_1$ of $A$. Similarly,
suppose that if $y'$ lies in $X_1$ and $A'$ is non-empty, then $y'$ 
is at least as large as the smallest element of $A'$. Then
$$c(A',y')(t(A,y)) = \cases{
1 & if $A = A'$ and $y = y'$ \cr
0 & otherwise.}$$
}

\noindent{\sl Proof.} Let us first consider the case where $A = \emptyset$.
Then $y \in X_2$, by assumption, and $t(A,y)$ is a single edge
$e$ labelled $y$ based at the basepoint of $\tilde L$. The evaluation
of $c(A',y')$ on $e$ is 
$$c(A',y')(e) = \left ( \prod_{j \in A'} c_j(i(e)) \right ) c_{y'}(e).$$
For this to be non-zero, we must have $y = y'$ and
$A' = \emptyset$. In this case, the evaluation is 1, as
required.

Let us now suppose that $A \not= \emptyset$ and hence that
$A' \not= \emptyset$. The loop $t(A,y)$ then consists of
two edges labelled $y$ and two edges labelled $y_1$.
Now, the evaluation of $c(A',y')$ on each such edge $e$
is again
$$c(A',y')(e) = \left ( \prod_{j \in A'} c_j(i(e)) \right ) c_{y'}(e).$$
This is zero if $c_{y'}(e) = 0$. Thus, the evaluation of
$t(A,y)$ is zero unless $y' = y$ or $y' = y_1$.
If neither of these equalities holds, the lemma is true.

\noindent {\sl Case 1.} $y' = y$.

Then, only the edges labelled $y$ contribute to the
evaluation of $t(A,y)$. Their initial vertices are
$v$ and $vy_1^{-1}$. Thus, the total evaluation is
$$\prod_{j \in A'} c_j(v) - \prod_{j \in A'} c_j(vy_1^{-1}).$$
The first term is 1 if and only if $A' \subseteq A$.
Since we are assuming $|A'| \geq |A|$, this happens
if and only if $A = A'$. Otherwise, the first term
is zero. The latter term is always zero, since 
to be non-zero, it would have to be the case that
$A' \subseteq A - \{ y_1 \}$. This proves the lemma in
the first case.

\noindent {\sl Case 2.} $y' = y_1$.

In this case, only the edges labelled $y_1$ contribute
to the evaluation of $t(A,y)$. Their initial vertices are
$vy_1^{-1}$ and $vyy_1^{-1}$. Thus, the total evaluation is 
$$\prod_{j \in A'} c_j(vy_1^{-1}) - \prod_{j \in A'} c_j(vyy_1^{-1}).$$
The first term is zero since $A' \not\subseteq A - \{ y_1 \}$.
If the second term is non-zero, then $A' \subseteq A \cup \{y \} - \{ y_1 \}$.
But comparing the sizes of these sets, we see that this
must be an equality. Hence, the smallest element
of $A'$ is strictly bigger than $y_1$, which equals $y'$. (Recall that
$y_1$ is the smallest element of $A$, and $y$ is larger than $y_1$.) We therefore
deduce that $y'$ is strictly smaller than every element of $A'$,
which is contrary to hypothesis. Thus, in Case 2, 
the evaluation of $t(A,y)$ is zero. $\square$

\noindent {\bf Lemma 3.13.} {\sl Let $p = 2$. Let $A'$ be a subset
of $X_1$ with size at most $\ell$, and let $y'$ be an element of
$X_1 \cup X_2$. Let $A$ be a subset of $X_1$ with size at
most $|A'|$, and let $y$ be an element of $A$. Then
$$c(A',y')(t(A,y)) = \cases{
1 & if $A = A'$ and $y = y'$ \cr
0 & otherwise.}$$}

\noindent {\sl Proof.} The test loop $t(A,y)$ has two edges
labelled $y$, with initial vertices $v$ and $vy$.
So, its evaluation under $c(A', y')$ is
$$\left ( \prod_{j \in A'} c_j(v) + \prod_{j \in A'} c_j(vy) \right ) c_{y'}(y).$$
This is zero unless $y' = y$. The first term in the brackets is
zero unless $A' \subseteq A$, which happens if and only if $A' = A$.
The second term is always zero, since $A' \not\subseteq A - \{ y \}$. $\square$

The reason why Theorem 1.6 treats odd primes in a different
way from the prime 2 arises in the above lemma. For odd $p$,
one may also define test loops $t(A,y)$ when $y \in A$:
this is a loop based at a suitable vertex $v$ (depending on $A$)
running along the edges $y^p$. However, the evaluation of
$c(A,y)$ on this $t(A,y)$ is $\sum_{r=0}^{p-1} r$, which is
zero modulo $p$, when $p$ is odd.

\noindent {\sl Proof of Proposition 3.10.}
We pick a total order on the test loops $T$, subject to 
the condition that if $|A| < |A'|$, then
$t(A,y) < t(A', y')$, whenever these loops are defined.
For each test loop $t' = t(A',y')$, let $z_{t'}$ be the cochain $c(A',y')$
in $U^1_\ell$. Then, by Lemmas 3.12 and 3.13,
$z_{t'}(t') = 1$ and $z_{t'}(t) = 0$ for all $t < t'$. 
This proves the proposition. $\square$

\noindent {\bf Proposition 3.14.} {\sl The dimension of
$C^1_\ell/(B^1 \cap C^1_\ell)$ is at least
$$\sum_{r=2}^{\ell+1} \left ( n \atop r \right )(r-1) + 
(b_1(G; {\Bbb F}_p)-n) \sum_{r=0}^{\ell}
\left ( n \atop r \right )  - b_2(G; {\Bbb F}_p) \sum_{r=0}^{\ell-1}
\left ( n \atop r \right )$$
if $p$ is odd, and at least
$$b_1(G; {\Bbb F}_2)\sum_{r=0}^\ell \left ( n \atop r \right ) 
- \sum_{r=1}^{\ell+1} \left ( n \atop r \right )  - 
b_2(G; {\Bbb F}_2) \sum_{r=0}^{\ell-1}
\left ( n \atop r \right )$$
if $p=2$.}

\noindent {\sl Proof.} We have the isomorphisms
$${C^1_\ell \over C^1_\ell \cap B^1} =
{C^1_\ell \over C^1_\ell \cap \psi(B^1)} =
{\psi({\rm ker}(\phi_\ell)) \over \psi({\rm ker}(\phi_\ell)) \cap \psi(B^1)}
\cong {{\rm ker}(\phi_\ell) \over {\rm ker}(\phi_\ell) \cap B^1}.$$
The first equality is a consequence of Lemma 3.7.
The second is just the definition of $C^1_\ell$.
The final isomorphism is a consequence of the
fact that $\psi$ is injective. Thus,
$$\eqalign{\dim(C^1_\ell / (C^1_\ell \cap B^1)) 
&\geq {\rm dim}({\rm ker}(\phi_\ell)) - {\rm dim}(U^1_\ell \cap B^1) \cr
&\geq {\rm dim}(U^1_\ell) - {\rm dim}(U^1_\ell \cap B^1)
- {\rm dim}(V_\ell).}$$
The proposition now follows from Proposition 3.10,
which gives a lower bound on the dimension of $U^1_\ell/(U^1_\ell 
\cap B^1)$, and the formula for the dimension of $V_\ell$.
$\square$

Theorem 1.6 immediately follows from this proposition and Theorem 3.8, since
$C^1_\ell/(B^1 \cap C^1_\ell)$ is a subspace of $H^1(\tilde L; {\Bbb F}_p)$.

\vskip 18pt
\centerline{\caps 4. Relationship with the lower central $p$-series}
\vskip 6pt

The proof of Theorem 1.6 was fairly formal. In this section, we
aim to explain it in terms that are possibly more familiar.

It is instructive to consider the case where  $\ell = 1$, $p=2$ and
$n = b_1(G; {\Bbb F}_2)$ in Theorem 1.6.
This forces $K$ to be $\gamma_2(G)$, the second
term in the lower-central 2-series for $G$.
Theorem 1.6 gives the inequality
$$b_1(\gamma_2(G); {\Bbb F}_2) \geq \left ( b_1(G; {\Bbb F}_2) \atop 2 \right )
+ b_1(G; {\Bbb F}_2) - b_2(G; {\Bbb F}_2).$$
This lower bound on $b_1(\gamma_2(G); {\Bbb F}_2)$ was already
known. Indeed, the following result was proved by
Shalen and Wagreich (see Lemma 1.3 of [6]).

\noindent {\bf Theorem 4.1.} {\sl Let $G$ be a group and let $p$ be a
prime. Suppose that $b_1(G; {\Bbb F}_p)$ and $b_2(G; {\Bbb F}_p)$
are finite. Then
$$b_1(\gamma_2(G); {\Bbb F}_p) \geq \left ( b_1(G; {\Bbb F}_p) \atop 2 \right )
+ b_1(G; {\Bbb F}_p) - b_2(G; {\Bbb F}_p).$$}

They proved this using the following
exact sequence of Stallings [7]:
$$H_2(G; {\Bbb F}_p) \rightarrow H_2(G/\gamma_2(G); {\Bbb F}_p) \rightarrow
\gamma_2(G)/\gamma_3(G) \rightarrow 0.$$
Now, $G/\gamma_2(G)$ is an elementary abelian $p$-group with rank
$b_1(G; {\Bbb F}_p)$. Its second homology is well-known to have
rank
$$\left ( b_1(G; {\Bbb F}_p) \atop 2 \right ) + b_1(G; {\Bbb F}_p),$$
via the K\"unneth formula. Thus, exactness of the sequence
gives that $\gamma_2(G)/\gamma_3(G)$ has dimension at least
$$\left ( b_1(G; {\Bbb F}_p) \atop 2 \right )
+ b_1(G; {\Bbb F}_p) - b_2(G; {\Bbb F}_p).$$
Since $\gamma_2(G)/\gamma_3(G)$ is a quotient of
$H_1(\gamma_2(G); {\Bbb F}_p)$, we deduce the required
lower bound on $b_1(\gamma_2(G); {\Bbb F}_p)$.

Now, $\gamma_2(G)/\gamma_3(G)$ is an elementary abelian $p$-group,
and hence it is isomorphic to $(\gamma_2(G)/\gamma_3(G))^\ast$,
which is
${\rm Hom}(\gamma_2(G)/\gamma_3(G); {\Bbb F}_p)$.
It is often useful to work with this latter group.
Any homomorphism $\gamma_2(G)/\gamma_3(G) \rightarrow {\Bbb F}_p$
arises from a homomorphism $\gamma_2(G) \rightarrow {\Bbb F}_p$
that is trivial on $\gamma_3(G)$. Conversely, any such
homomorphism gives an element of $(\gamma_2(G)/\gamma_3(G))^\ast$.
Thus, one can consider $(\gamma_2(G)/\gamma_3(G))^\ast$
to be a subgroup of the set of all homomorphisms from
$\gamma_2(G)$ to ${\Bbb F}_p$. This is just $H^1(\gamma_2(G); {\Bbb F}_p)$.
Now, $H^1(\gamma_2(G); {\Bbb F}_p)$ is isomorphic to $H^1(\tilde L; {\Bbb F}_p)$,
where $\tilde L$ is the 2-complex from Section 3. The Stallings
exact sequence gives a lower bound on the dimension of
$(\gamma_2(G)/\gamma_3(G))^\ast$. Each element in here gives
an element of $H^1(\tilde L; {\Bbb F}_p)$, which is represented
by a 1-cocycle on $\tilde L$. What are these cocycles? When $p=2$, they
are precisely $C^1_\ell$ for $\ell = 1$. 

Thus, the level $\ell = 1$ is the topological analogue of 
$\gamma_2(G)/\gamma_3(G)$. As one might expect, higher values
of $\ell$ do indeed correspond to sections further down
the lower central $p$-series of $G$. Specifically, one
can consider the covering space $\tilde L_{\ell+1}$ corresponding
to the subgroup $\gamma_{\ell +1}(G)$ of $G$. It is
possible to construct explicit 1-cocycles on $\tilde L_{\ell +1}$
representing certain elements of $(\gamma_{\ell+1}(G)/\gamma_{\ell +2}(G))^\ast$.
These cocycles turn out to be invariant under the action of
the covering group $\gamma_2(G)/\gamma_{\ell + 1}(G)$ and so
descend to cocycles on $\tilde L$. These lie in
$C^1_\ell$.

Thus, Theorem 1.6 arose from an attempt to understand the Stallings
exact sequence topologically, and to explore its possible analogues
further down the lower central $p$-series. So the appearance at
various stages in Sections 2 and 3 of $\gamma_m(F)$ was more
than just a technical device. The lower central $p$-series has a 
crucial r\^ole in the interpretation of Theorem 1.6.

\vskip 18pt
\centerline {\caps 5. Normal subgroups with large homology}
\vskip 6pt

Most of the theorems in this paper make the hypothesis that the group $G$ contains
finite index subgroups $G_1$ where $b_1(G_1; {\Bbb F}_p)$ is arbitrarily
large. In this section, we show that we may assume that these
subgroups $G_1$ are, in addition, normal. We will need to make a 
hypothesis about $G$ that is much weaker than the $b_2 - b_1$ 
and $b_2/b_1$ conditions.

This section is not in fact required for most of the results in this paper.
It is necessary only to prove that the subgroups
$G_i$ in Theorem 1.7 with fast homology growth are normal in $G$. If one is
content with the weaker conclusion that they are just normal in
$G_1$, then this section could be omitted entirely. From this,
there is an easy argument which gives that we may take each $G_i$
to be subnormal in $G$. This would be sufficient to deduce Theorem 1.1.
However, we prefer to pursue the strongest possible conclusion for
Theorem 1.7: that each $G_i$ is normal in $G$. For this, it appears
that more technology is required: we need some results about $p$-adic analytic
pro-$p$ groups.

Our main result in this section is the following.

\noindent {\bf Theorem 5.1.} {\sl Let $G$ be a finitely generated
group and let $p$ be a prime. Suppose that, for some finite index
subgroup $G_1$ of $G$, $b_1(G_1; {\Bbb F}_p) > 1$ and $b_2(G_1;{\Bbb F}_p) <
b_1(G_1;{\Bbb F}_p)^2/4$. Then
$$\sup \{ b_1(G_i; {\Bbb F}_p): G_i \hbox{ is a finite index
normal subgroup of } G \} = \infty.$$}

We will need to quote two facts about $p$-adic analytic
pro-$p$ groups. For the following, see Interlude D in [1].

\noindent {\bf Theorem 5.2.} {\sl Let $G$ be a finitely generated group
and let $p$ be a prime.
Suppose that  $b_1(G; {\Bbb F}_p) > 1$ and
$b_2(G;{\Bbb F}_p) < b_1(G;{\Bbb F}_p)^2/4$.
Then the pro-$p$ completion of $G$, denoted $\hat G_{(p)}$, is
not $p$-adic analytic.}

A proof of the following can be found in [4].

\noindent {\bf Theorem 5.3.} {\sl Let $G$ be a finitely generated
group and let $p$ be a prime. Then the following are equivalent:
\item{1.} $\hat G_{(p)}$ is $p$-adic analytic;
\item{2.} the supremum of $b_1(G_i; {\Bbb F}_p)$, as $G_i$ ranges over all
characteristic subgroups of $G$ with index a power of $p$, is finite;
\item{3.} the supremum of $b_1(G_i; {\Bbb F}_p)$, as $G_i$ ranges over all
normal subgroups of $G$ with index a power of $p$, is finite.

}

\noindent {\sl Proof of Theorem 5.1.}
By hypothesis, there is a finite index subgroup $G_1$ of $G$ 
such that $b_1(G_1; {\Bbb F}_p) > 1$ and
$b_2(G_1;{\Bbb F}_p) < b_1(G_1;{\Bbb F}_p)^2/4$.
So, by Theorem 5.2, the pro-$p$ completion of
$G_1$ is not $p$-adic analytic. Therefore, by Theorem 5.3,
$G_1$ contains a sequence of characteristic subgroups
$G_i$, each with index a power of $p$, such that
$b_1(G_i; {\Bbb F}_p)$ tends to infinity.

Let $K$ be a finite index normal subgroup of $G$ that lies in
$G_1$. Let $K_i$ be the intersection of $K$ and $G_i$.
Then,
$$b_1(K_i; {\Bbb F}_p) \geq b_1(G_i; {\Bbb F}_p) - b_1(G_i/K_i; {\Bbb F}_p).$$
But $G_i/K_i$ is isomorphic to $G_iK/K$, which is
a subgroup of $G_1/K$. Since this finite group has only
finitely many subgroups, $b_1(G_i/K_i; {\Bbb F}_p)$ is uniformly bounded
above. Hence, $b_1(K_i; {\Bbb F}_p)$ tends to infinity. Now, $G_i$ is
normal in $G_1$ and has index a power of $p$. So, $K_i$ is normal in $K$
and has index a power of $p$. Hence,
by Theorem 5.3, $K$ has a sequence of finite index
characteristic subgroups $L_i$ such that $b_1(L_i; {\Bbb F}_p)$
tends to infinity. Since these are characteristic in $K$,
which is normal in $G$, they are therefore normal in $G$.
These are the required subgroups of $G$. $\square$

\vskip 18pt
\centerline {\caps 6. Homology growth and subgroup growth}
\vskip 6pt

In this section, we use the homological lower bounds of Theorem 1.6
to deduce Theorem 1.7. We are assuming that $G$ is a finitely generated group
having the $b_2 - b_1$ property with respect to the prime $2$ and that
$$\sup \{ b_1(G_i; {\Bbb F}_2) : G_i \hbox{ is a finite index subgroup of } G \} = \infty.$$
Pick $G_1$ where $b_1(G_1; {\Bbb F}_2)$ is bigger than 1 and large enough so that
$b_2(G_1;{\Bbb F}_p) < b_1(G_1;{\Bbb F}_p)^2/4$.
Theorem 5.1 states that
$$\sup \{ b_1(G_i; {\Bbb F}_p): G_i \hbox{ is a finite index
normal subgroup of } G \} = \infty.$$
Thus, we may assume not only that $b_1(G_1; {\Bbb F}_p)$ is large
but also that $G_1$ is normal in $G$.

Define subgroups $G_i$ of $G$ recursively, by
setting $G_{i+1} = [G_i,G_i](G_i)^2$. Thus,
$G_i$ is just the derived 2-series for $G_1$.
Let $x_i = b_1(G_i; {\Bbb F}_2)$.
According to Theorem 1.6, setting $\ell = \lfloor x_i/2 \rfloor$
and $n = x_i$, we have
$$
x_{i+1} \geq x_i \sum_{r=0}^{\lfloor {x_i}/2 \rfloor} \left ( x_i \atop r \right ) 
- \sum_{r =1}^{\lfloor x_i/2 \rfloor +1} \left ( {x_i \atop r} \right ) - 
b_2(G_i; {\Bbb F}_2) \sum_{r=0}^{\lfloor x_i/2 \rfloor - 1}
\left ( x_i \atop r \right ).$$
Now, the middle summation is bounded below by $-2^{x_i}$.
The third summation can be compared with all but the
highest term in the first summation. Thus,
$$x_{i+1} \geq x_i \left ( x_i \atop \lfloor x_i/2 \rfloor \right ) 
- 2^{x_i} \max \{ 1, 1 + b_2(G_i; {\Bbb F}_2) - x_i \}.\eqno{(1)}$$

\noindent {\sl Claim 1.} Let $\lambda$ be any positive real number less
than $\sqrt{2/\pi}$. Then, provided $x_1$ is sufficiently big,
$$x_{i+1} \geq \lambda 2^{x_i} \sqrt{x_i},$$
for all $i \geq 1$.

According to Stirling's formula,
$$x! \sim \sqrt{2 \pi x} \ x^x e^{-x},$$
as $x \rightarrow \infty$. So,
$$\left( {x \atop \lfloor x/2 \rfloor} \right ) 
\sim {\sqrt{2 \pi x} x^x e^{-x} 
\over 2 \pi (x/2)^{x+1} e^{-x}}
= \sqrt{2 \over \pi} {2^x \over \sqrt x}.$$
Thus, the first term in (1) is at least
$$\lambda 2^{x_i} \sqrt{x_i},$$
when $x_i$ is sufficiently large.
Now, by the $b_2 - b_1$ condition, there is a universal upper
bound on $b_2(G_i; {\Bbb F}_2) - x_i$. Thus, the first term of (1)
dominates, and the claim is proved. In fact, to prove this claim,
one does not need the full strength of the $b_2 - b_1$ condition.
It suffices to assume that, when $x_i$ is sufficiently large,
$b_2(G_i; {\Bbb F}_2) - x_i$ is negative or small compared with
$\sqrt{x_i}$. 

Note that, by the claim, if we pick $b_1(G_1; {\Bbb F}_2)$
to be sufficiently big, then $x_i$ is a strictly increasing
function.

Set $$\sigma_i = \sum_{j=1}^i x_j.$$ 

\noindent {\sl Claim 2.} Provided $x_1$ is sufficiently big, then
for all $i \geq 1$,
$$2^{\sigma_i} \leq \lambda 2^{x_i} x_i (\log x_i)^{2/3}.$$
We prove this by induction on $i$.  It is clear for $i =1$.
For the inductive step, note that
$$\eqalign{
2^{\sigma_{i+1}}
&= 2^{x_{i+1}} 2^{\sigma_i} \cr
&\leq 2^{x_{i+1}} \lambda 2^{x_i} x_i (\log x_i)^{2/3} \cr
&\leq 2^{x_{i+1}} x_{i+1} x_i^{1/2} (\log x_i)^{2/3} \cr
&\leq 2^{x_{i+1}} x_{i+1} (\log x_{i+1})^{2/3},}$$
where the second inequality is consequence of Claim 1 and
the final step follows from the fact that
$$(\log x_{i+1})^{2/3} \geq (x_i \log 2 + \log \sqrt{x_i} + \log \lambda)^{2/3} 
\geq x_i^{1/2} \log x_i,$$
provided $x_i$ is sufficiently big. This certainly holds if
$x_1$ is sufficiently large.

\noindent {\sl Claim 3.} As $i \rightarrow \infty$,
$$x_{i} = \Omega \left ( {2^{\sigma_{i-1}} \over \sigma_{i-1}^{1/2} \log \sigma_{i-1}} \right ).$$

This follows from
$$x_{i} \geq \lambda 2^{x_{i-1}} x_{i-1}^{1/2} \geq {2^{\sigma_{i-1}} \over x_{i-1}^{1/2} 
(\log x_{i-1})^{2/3}}
\geq {2^{\sigma_{i-1}} \over \sigma_{i-1}^{1/2} (\log \sigma_{i-1})^{2/3}} = 
\Omega \left ( {2^{\sigma_{i-1}} \over \sigma_{i-1}^{1/2} \log \sigma_{i-1}} \right ).$$
The second inequality is an application of Claim 2, and the
penultimate inequality comes from the fact that
$\sigma_{i-1} \geq x_{i-1}$.

The proof of Theorem 1.7 follows quickly from this claim.
This is because $x_i = b_1(G_i; {\Bbb F}_2)$, by definition,
and 
$$[G_1:G_i] = \prod_{j=1}^{i-1} [G_j:G_{j+1}] = 2^{\sum_{j=1}^{i-1} x_j}
= 2^{\sigma_{i-1}}.$$
This completes the proof of Theorem 1.7. We saw in the introduction
that Theorem 1.1 is a rapid consequence of Theorem 1.7.

It should now be apparent why Theorems 1.1 and 1.7 work only
when $p =2$. This is {\sl not} due to the two different
lower bounds in Theorem 1.6. Instead, it is a consequence of
the behaviour of binomial coefficients. If one wants to establish
fast homology growth for the sequence of subgroups
$G_i$, one needs to know that the ratio $b_1(G_{i+1}; {\Bbb F}_p)/
b_1(G_i; {\Bbb F}_p)$ is comparable with the index
$[G_i:G_{i+1}]$. Now, the latter is 
$$p^{b_1(G_i; {\Bbb F}_p)}.$$
A lower bound on the former comes from Theorem 1.6. When
$b_2(G_i; {\Bbb F}_p) = b_1(G_i; {\Bbb F}_p)$, say,
the formula in Theorem 1.6 is at most
$$b_1(G_i; {\Bbb F}_p) \left ( b_1(G_i; {\Bbb F}_p) \atop \ell \right ).$$
So, no matter what value of $\ell$ is chosen,
the lower bound on $b_1(G_{i+1}; {\Bbb F}_p)/
b_1(G_i; {\Bbb F}_p)$ that we obtain is at most 
$$2^{b_1(G_i; {\Bbb F}_p)}.$$
Thus, the only situation in which we can prove that the growth in homology is
comparable to the growth in the subgroups' index is when
$p = 2$. If one wanted to prove similar results when $p$ is
odd, a significantly stronger version of Theorem 1.6 would be
required.

It also clear that, in practice, one does need the
full hypotheses of Theorem 1.1 to deduce its conclusion.
Firstly, as noted in the proof of Claim 1, one may weaken the $b_2 - b_1$ condition
to the hypothesis that, for any sequence of finite index subgroups
$G_i$ of $G$ such that $b_1(G_i; {\Bbb F}_2) \rightarrow \infty$,
$$\limsup_i {b_2(G_i; {\Bbb F}_2) - b_1(G_i; {\Bbb F}_2) \over \sqrt{b_1(G_i; {\Bbb F}_2)}} \leq 0.$$

Secondly, one does not need to {\sl assume} that
$$\sup \{ b_1(G_i; {\Bbb F}_2) : G_i \hbox{ is a finite index subgroup of } G \} = \infty.$$
This is in fact a consequence of the $b_2 - b_1$ condition, together
with the fact that $b_1(G_1; {\Bbb F}_2)$ is `sufficiently large',
for some finite index subgroup $G_1$ of $G$. For, suppose that one
has an upper bound ($k$, say) on $b_2(G_i; {\Bbb F}_2) - b_1(G_i; {\Bbb F}_2)$
for all finite index subgroups $G_i$. One can then apply Theorem 1.6
with $\ell = 1$ (or Theorem 4.1)
to a finite index subgroup $G_1$ of $G$ to give
$$b_1(\gamma_2(G_1); {\Bbb F}_2) \geq \left ( b_1(G_1; {\Bbb F}_2) \atop 2 \right )
+ b_1(G_1; {\Bbb F}_2) - b_2(G_1; {\Bbb F}_2) \geq \left ( b_1(G_1; {\Bbb F}_2) \atop 2 \right )
- k.$$
Suppose $b_1(G_1; {\Bbb F}_2)$ is big enough (depending on $k$)
to ensure that $b_1(\gamma_2(G_1); {\Bbb F}_2) > b_1(G_1; {\Bbb F}_2)$.
We may then repeat to find a lower bound on $b_1(\gamma_2(\gamma_2(G_1)); {\Bbb F}_2)$,
and so on. Thus, we obtain a sequence of finite index subgroups
$\{ G_i \}$ such that $b_1(G_i; {\Bbb F}_2)$ tends to infinity.
Moreover, this is the derived $2$-series of some finite index
subgroup $G_1$ of $G$. We will see this approach put into practice in
Section 8 with the proofs of Theorems 1.4 and 1.5.

This then leads to a third way that
the hypotheses of Theorem 1.1 may be weakened. One does not
need to assume a uniform upper bound on $b_2(G_i; {\Bbb F}_2) - b_1(G_i; {\Bbb F}_2)$
over all finite index subgroups $G_i$. Once one has a finite index
subgroup $G_1$ of $G$ such that $b_1(G_i; {\Bbb F}_2)$ tends to
infinity for the derived 2-series $\{ G_i \}$ for $G_1$, then one
only needs to assume an upper bound on $b_2(G_i; {\Bbb F}_2) - b_1(G_i; {\Bbb F}_2)$
for this collection of finite index subgroups.

\vfill\eject
\centerline{\caps 7. The $b_2/b_1$ condition}
\vskip 6pt

Let $G$ be a group satisfying the
$b_2/b_1$ condition with respect to the prime $p$. 
Let $m$ be the supremum of
$$\left \lceil {b_2(G_i; {\Bbb F}_p) \over b_1(G_i; {\Bbb F}_p)} \right \rceil$$
as $G_i$ ranges over all finite index subgroups 
such that $b_1(G_i; {\Bbb F}_p) > 0$. Then for 
all such $G_i$, $b_2(G_i; {\Bbb F}_p) \leq m \ b_1(G_i; {\Bbb F}_p)$.
Suppose also that
$$\sup \{ b_1(G_i; {\Bbb F}_p) : G_i \hbox{ is a finite index subgroup of } G \} = \infty.$$
We will prove that there is a sequence
of finite index subgroups $G \triangleright G_1 \triangleright G_2 \triangleright \dots$
with the following property. For all $\epsilon > 0$,
$$b_1(G_i; {\Bbb F}_p) \geq [G:G_i]^{(1/(m+p-1)) - \epsilon}$$
for all sufficiently large $i$.
This will prove Theorems 1.8 and 1.2. It also
provides a lower bound on the constant $k$ in these
theorems.

By Theorem 5.1, there are finite index normal subgroups $G_1$ of $G$
where $b_1(G_1; {\Bbb F}_p)$ is arbitrarily large.
Set $n = m + p -1$ and set $\ell = 1$. Define a sequence of subgroups
$\{ G_i \}$ of $G$, as follows. For $i \geq 1$, let
$G_{i+1}$ be a normal subgroup of $G_i$
such that $G_i/G_{i+1}$ is an elementary abelian
$p$-group of rank $n$. Applying Theorem 1.6, obtain
the inequality
$$\eqalign{
b_1(G_{i+1}; {\Bbb F}_p) & \geq  \left ( n \atop 2 \right ) + 
(b_1(G_i; {\Bbb F}_p) - n)(n+1) - b_2(G_i; {\Bbb F}_p) \cr
&\geq b_1(G_i; {\Bbb F}_p) p - n^2/2 - 3n/2.}$$
Now, when $b_1(G_i; {\Bbb F}_p)$ is sufficiently
large, the $n^2/2$ and $3n/2$ terms are negligible compared
with the first term. In particular, $b_1(G_{i+1}; {\Bbb F}_p)$
is strictly greater than $b_1(G_i; {\Bbb F}_p)$.
Thus, we may repeat this indefinitely, obtaining
a sequence of finite index subgroups $G_1 \triangleright G_2 \triangleright \dots$
such that $b_1(G_i; {\Bbb F}_p) \rightarrow \infty$.
Thus,
$$\liminf_i {b_1(G_{i+1}; {\Bbb F}_p) \over b_1(G_i; {\Bbb F}_p)} \geq p.$$
But, $[G_i:G_{i+1}] = p^{(m+p-1)}$.
Thus, 
$$\liminf_i {(\log_p b_1(G_{i+1}; {\Bbb F}_p) - \log_p b_1(G_i; {\Bbb F}_p)) 
\over \log_p [G_i:G_{i+1}]}
\geq {1 \over m+p-1}.$$
Thus,
$$\eqalign{
\liminf_i {\log_p b_1(G_i; {\Bbb F}_p) \over \log_p [G:G_i]} 
&\geq \liminf_i {\sum_{j=1}^{i-1} (\log_p b_1(G_{j+1}; {\Bbb F}_p) - \log_p b_1(G_j; {\Bbb F}_p))
\over \log_p [G:G_1] + \sum_{j=1}^{i-1} \log_p [G_j:G_{j+1}]} \cr
&\geq {1 \over m+p-1}.}$$
This implies that, for all $\epsilon > 0$, 
$$b_1(G_i; {\Bbb F}_p) \geq {[G:G_i]^{(1/(m+p-1)) - \epsilon}},$$
for all sufficiently large $i$. $\square$

\vskip 18pt
\centerline{\caps 8. Examples}
\vskip 6pt

In this section, we show that the $b_2 - b_1$ and $b_2/b_1$ conditions
hold for several different classes of groups. The aim is to establish
Theorems 1.4 and 1.5.

\noindent {\bf Proposition 8.1.} {\sl Let $G$ be the fundamental group
of a compact 3-manifold $M$ that is either closed or orientable,
and let $p$ be a prime. Then $b_2(G_i; {\Bbb F}_p) \leq 
b_1(G_i; {\Bbb F}_p)$ for any finite index subgroup $G_i$ of $G$.}

\noindent {\sl Proof.} Any finite index subgroup $G_i$ of $G$ is the
fundamental group of a compact 3-manifold $\tilde M$ that is either
closed or orientable. By attaching 3-balls to $\tilde M$ if necessary, we
may assume that it has no 2-sphere boundary components. Thus,
its Euler characteristic is non-positive, by Poincar\'e duality. So,
$$\eqalign{
0 &\geq \chi(\tilde M) = b_0(\tilde M; {\Bbb F}_p) 
- b_1(\tilde M; {\Bbb F}_p) + b_2(\tilde M; {\Bbb F}_p) - b_3(\tilde M; {\Bbb F}_p) \cr
&\geq - b_1(\tilde M; {\Bbb F}_p) + b_2(\tilde M; {\Bbb F}_p).}$$
The proposition is therefore a consequence of the following
easy fact. $\square$

\noindent {\bf Lemma 8.2.} {\sl Let $X$ be a topological space with
fundamental group $G$. Then
$$\eqalign{
b_1(G; {\Bbb F}_p) &= b_1(X; {\Bbb F}_p) \cr
b_2(G; {\Bbb F}_p) &\leq b_2(X; {\Bbb F}_p).}$$
}

\noindent {\sl Proof.} An Eilenberg-Maclane space $K(G,1)$ can
be constructed from $X$ by attaching cells with dimensions
at least 3. This does not affect $b_1$ and it does not increase
$b_2$. $\square$

Any lattice in ${\rm PSL}(2, {\Bbb C})$ has a finite index subgroup
that is the fundamental group of a compact orientable 3-manifold. So,
Proposition 8.1, Theorem 1.3 and Theorem 1.1
deal with case (1) of Theorem 1.4. The following deals with case (2).

\noindent{\bf Proposition 8.3.} {\sl Let $G$ and $p$ be as in Proposition 8.1.
Suppose that $b_1(G_i; {\Bbb F}_p) > 3$ for some finite index
subgroup $G_i$ of $G$. Then
$$\sup \{ b_1(G_i; {\Bbb F}_p) : G_i \hbox{ is a finite index subgroup of } G \} = \infty.$$
}

\noindent {\sl Proof.} By Theorem 4.1 and Proposition 8.1,
$$b_1(\gamma_2(G_i); {\Bbb F}_p) \geq \left ( {b_1(G_i; {\Bbb F}_p) \atop 2} \right )
+ b_1(G_i; {\Bbb F}_p) - b_2(G_i; {\Bbb F}_p) \geq
\left ( {b_1(G_i; {\Bbb F}_p) \atop 2} \right ).$$
This is strictly greater than $b_1(G_i; {\Bbb F}_p)$, since
$b_1(G_i; {\Bbb F}_p) > 3$. Repeating for $\gamma_2(\gamma_2(G_i))$,
and so on, we obtain a sequence $\{ G_j \}$ of finite index subgroups of
$G$, such that $b_1(G_j; {\Bbb F}_p)$ tends to infinity. $\square$

Note that the hypothesis that $b_1(G_i; {\Bbb F}_p) > 3$ is
necessary here. For example, when $G$ is the fundamental group
of the 3-torus, then any finite index subgroup of $G$ is
isomorphic to ${\Bbb Z} \times {\Bbb Z} \times {\Bbb Z}$.
Note also that here, $G$ has polynomial subgroup growth.

We now consider groups with deficiency at least 1. By definition,
these are groups with a finite presentation $\langle X | R \rangle$
where $|X| - |R| = 1$. Note that, by the Reidemeister-Schreier
process, any finite index subgroup of a group with deficiency
at least 1 also has deficiency at least 1.

\noindent{\bf Proposition 8.4.} {\sl Let $G$ be a group with
deficiency at least 1. Then 
$$b_2(G; {\Bbb F}_p) - b_1(G; {\Bbb F}_p) \leq -1.$$}

\noindent {\sl Proof.} Let $L$ be the 2-complex arising
from the presentation $\langle X|R \rangle$
for $G$ where $|X| - |R| = 1$. Then
$$b_2(L; {\Bbb F}_p) - b_1(L; {\Bbb F}_p) + 1 = 
\chi(L) = |R| - |X| + 1 = 0.$$
Now apply Lemma 8.2. $\square$

\noindent {\bf Corollary 8.5.} {\sl Any group with 
deficiency at least 1 satisfies the $b_2 - b_1$ condition.}

\noindent {\bf Proposition 8.6.} {\sl Let $G$ be a group
with deficiency at least 1.
Suppose that $b_1(G_i; {\Bbb F}_p) > 2$ for some finite index
subgroup $G_i$ of $G$. Then
$$\sup \{ b_1(G_i; {\Bbb F}_p) : G_i \hbox{ is a finite index subgroup of } G \} = \infty.$$
}

\noindent {\sl Proof.} This is essentially the same as
the proof of Proposition 8.3, but we use Proposition 8.4 rather
than Proposition 8.1. $\square$

This deals with case (3) of Theorem 1.4. For case (4), note that
any (finitely generated free)-by-cyclic group $F \rtimes {\Bbb Z}$ has deficiency
at least 1. Thus, by Corollary 8.5, it satisfies the $b_2 - b_1$
condition. We also have the following, which completes
the analysis of case (4) of Theorem 1.4.

\noindent {\bf Proposition 8.7.} {\sl Let $G$ be a (finitely generated
free non-abelian)-by-cyclic group. Then, for any prime $p$,
$$\sup \{ b_1(G_i; {\Bbb F}_p) : G_i \hbox{ is a finite index subgroup of } G \} = \infty.$$
}

\noindent {\sl Proof.} The group $G$ is a semi-direct product $F \rtimes_\phi {\Bbb Z}$,
determined by an automorphism $\phi$ of the finitely generated free non-abelian group $F$.
Let $K$ be any finite index characteristic subgroup of $F$. This is 
preserved by $\phi$. It is clear that $K \rtimes_{\phi|K} {\Bbb Z}$
is a finite index subgroup of $F \rtimes_\phi {\Bbb Z}$.
In this way, we may assume that $F$ has arbitrarily large
rank. Now, $\phi$ induces an automorphism $\phi_\ast \colon H_1(F; {\Bbb F}_p)
\rightarrow H_1(F; {\Bbb F}_p)$. Since $H_1(F; {\Bbb F}_p)$
has finite order, $\phi_\ast^n$ is the identity for some
positive integer $n$. The kernel of the map $F \rtimes_{\phi} {\Bbb Z}
\rightarrow {\Bbb Z} \rightarrow ({\Bbb Z}/n{\Bbb Z})$ is
isomorphic to $F \rtimes_{\phi^n} {\Bbb Z}$. Note that
$b_1(F \rtimes_{\phi^n} {\Bbb Z}; {\Bbb F}_p) = b_1(F; {\Bbb F}_p) + 1$. Since
$b_1(F; {\Bbb F}_p)$ can be assumed to be arbitrarily large, the proposition
is proved. $\square$

We end with a large class of examples of $b_2/b_1$ groups.
These give case (3) of Theorem 1.5.

\noindent {\bf Proposition 8.8.} {\sl Let $G$ be the fundamental
group of a closed 4-manifold $M$ with non-positive Euler characteristic, 
and let $p$ be a prime. Suppose, in addition, that if $p$ is odd, then $M$ 
is orientable. Then $G$ satisfies the $b_2/b_1$ condition.}

\noindent {\sl Proof.} Corresponding to any finite index subgroup $G_i$
of $G$, there is a finite-sheeted covering space $\tilde M$ of $M$.
This also has non-positive Euler characteristic. When $M$ is orientable,
so is $\tilde M$. Now, Poincar\'e duality applied to $\tilde M$
gives that
$$\eqalign{
0 &\geq \chi(\tilde M) = b_0(\tilde M; {\Bbb F}_p)
- b_1(\tilde M; {\Bbb F}_p)
+ b_2(\tilde M; {\Bbb F}_p)
- b_3(\tilde M; {\Bbb F}_p)
+ b_4(\tilde M; {\Bbb F}_p) \cr
& = 2 - 2 b_1(\tilde M; {\Bbb F}_p) + b_2(\tilde M; {\Bbb F}_p).}$$
Thus, by Lemma 8.2,
$$b_2(G_i; {\Bbb F}_p) \leq b_2(\tilde M; {\Bbb F}_p)
\leq 2 b_1(G_i; {\Bbb F}_p) -2.$$
$\square$

\noindent {\bf Proposition 8.9.} {\sl Let $G$ and $p$ be as in Proposition 8.8.
Suppose that $b_1(G_i; {\Bbb F}_p) > 4$ for some finite index
subgroup $G_i$ of $G$. Then
$$\sup \{ b_1(G_i; {\Bbb F}_p) : G_i \hbox{ is a finite index subgroup of } G \} = \infty.$$
}

\noindent {\sl Proof.} This follows the same lines as the proof of
Proposition 8.3, using the inequality
$$b_2(G_i; {\Bbb F}_p) \leq 2 b_1(G_i; {\Bbb F}_p) - 2$$
that was established in the proof of Proposition 8.8. $\square$

Case (3) of Theorem 1.5 is proved by applying Propositions 8.8 and 8.9
and Theorem 1.2. When $M$ is orientable or $p=2$, the result follows
immediately. When $M$ is non-orientable and $p > 2$, we must first pass
to the orientable double cover $\tilde M$ of $M$. Let $K_i$ be
$\pi_1(\tilde M) \cap G_i$. Then
$$b_1(K_i; {\Bbb F}_p) \geq b_1(G_i; {\Bbb F}_p) - b_1(G_i/K_i; {\Bbb F}_p)
= b_1(G_i; {\Bbb F}_p) > 4.$$
Thus, $\pi_1(\tilde M)$ satisfies the hypotheses of Propositions 8.8 and 8.9.
Hence, Propositions 8.8 and 8.9 and Theorem 1.2 establish the theorem
in this case.

\vskip 18pt
\centerline{\caps References}
\vskip 6pt

\item{1.} {\caps J. Dixon, M. du Sautoy, A. Mann, D. Segal,}
{\sl Analytic pro-$p$ groups.} Cambridge Studies in Advanced Mathematics, 61. 
Cambridge University Press, Cambridge (1999).

\item{2.} {\caps P. Hall}, {\sl The Theory of Groups} (1959)

\item{3.} {\caps M. Lackenby}, {\sl Detecting large groups},
Preprint.

\item{4.} {\caps A. Lubotzky, A. Mann,}
{\sl Powerful $p$-groups II. $p$-adic analytic groups,}
J. Algebra, 105 (1987) 506--515.

\item{5.} {\caps A. Lubotzky, D. Segal},
{\sl Subgroup growth}. Progress in Mathematics, 212. 
Birkh\"auser Verlag (2003)

\item{6.} {\caps P. Shalen, P. Wagreich}, {\sl Growth rates, 
$Z\sb p$-homology, and volumes of hyperbolic $3$-manifolds. }
Trans. Amer. Math. Soc. 331 (1992) 895--917.

\item{7.} {\caps J. Stallings}, {\sl Homology and lower central series of
groups}, J. Algebra 2 (1965) 170--181.

\vskip 12pt
\+ Mathematical Institute, University of Oxford, \cr
\+ 24-29 St Giles', Oxford OX1 3LB, United Kingdom. \cr

\end